\documentclass[10pt,reqno]{amsart}
\usepackage{hannah}
\usepackage{color,graphicx}
\usepackage[latin1]{inputenc}
\DeclareMathOperator{\val}{val}

\newcommand {\brk}{\rule {1mm}{0mm}}
\newcommand {\df}[1]{\emph {#1}}

\newcommand {\calM}{{\mathcal M}}
\newcommand {\calP}{{\mathcal P}}
\DeclareMathOperator{\Log}{Log}

\begin{document}

\title [Tropical plane curves through points in general position]{The numbers
  of tropical plane curves through points in general position}
\author {Andreas Gathmann and Hannah Markwig}
\address {Andreas Gathmann, Fachbereich Mathematik, Technische Universität
  Kaiserslautern, Postfach 3049, 67653 Kaiserslautern, Germany}
\email {andreas@mathematik.uni-kl.de}
\address {Hannah Markwig, Fachbereich Mathematik, Technische Universität
  Kaiserslautern, Postfach 3049, 67653 Kaiserslautern, Germany}
\email {markwig@mathematik.uni-kl.de}
\thanks {\emph {2000 Mathematics Subject Classification:} Primary 14N10, 51M20,
  Secondary 14N35}
\thanks {The second author has been funded by the DFG grant Ga 636/2.}

\begin{abstract}
  We show that the number of tropical curves of given genus and degree through
  some given general points in the plane does not depend on the position of the
  points. In the case when the degree of the curves contains only primitive
  integral vectors this statement has been known for a while now, but the only
  known proof was indirect with the help of Mikhalkin's Correspondence Theorem
  that translates this question into the well-known fact that the numbers of
  \emph {complex} curves in a toric surface through some given points do not
  depend on the position of the points. This paper presents a direct proof
  entirely within tropical geometry that is in addition applicable to arbitrary
  degree of the curves.
\end{abstract}

\maketitle

%%%%%%%%%%%%%%%%%%%%%%%%%%%%%%%%%%%%%%%%%%%%%%%%%%%%%%%%%%%%%%%%%%%%%%%%%%%%%%%

\section{Introduction}

Tropical geometry recently attracted a lot of attention. One reason for this is
the possibility to relate \emph {complex} enumerative geometry to the
(hopefully simpler) \emph {tropical} enumerative geometry. For example,
Mikhalkin has proven the so-called ``Correspondence Theorem'' which asserts
that the numbers of complex curves (of given genus and homology class) in toric
surfaces through some given points are equal to the numbers of certain
plane tropical curves through the same number of given points \cite{Mi03}.
Furthermore, Mikhalkin gave a nice purely combinatorial way of computing these
numbers of tropical curves. Nishinou and Siebert were able to prove the same
for rational curves in (higher-dimensional) complete toric varieties
\cite{NS04}.

Of course the numbers of complex curves in a toric surface through some given
general points do not depend on the position of the points. It is therefore a
corollary of the Correspondence Theorem that the corresponding numbers of plane
tropical curves through some given general points cannot depend on the position
of the points either (see remark \ref{mik} for details). From a purely
tropical point of view this statement is far from being obvious however. It is
the main result of this paper to prove this statement within the framework of
tropical geometry. In addition, our result will be more general than what is
known so far since not all numbers of plane tropical curves have been related
to numbers of complex curves yet: for some numbers that should correspond to
relative Gromov-Witten invariants (i.e.\ curves with fixed local multiplicities
to a given divisor) there is no Correspondence Theorem yet (see remark
\ref{enden} for details).

Let us brief\/ly describe the tropical set-up. A plane tropical curve can
roughly be thought of as a weighted graph in the real plane whose edges are
(possibly unbounded) line segments of rational slope (for a precise definition
see section \ref{sec-planetrop}). These graphs must satisfy the so-called \emph
{balancing condition}: at each vertex the weighted sum of the primitive
integral vectors along the edges starting from this vertex must be zero. For
example, the following picture shows a tropical curve through 3 given points in
the plane:

\begin {center} \begin{picture}(0,0)%
\includegraphics{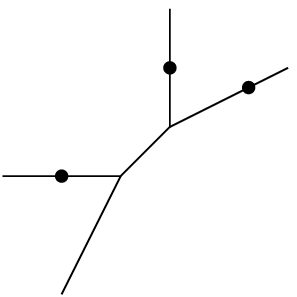}%
\end{picture}%
\setlength{\unitlength}{4144sp}%
\begingroup\makeatletter\ifx\SetFigFont\undefined
% extract first six characters in \fmtname
\def\x#1#2#3#4#5#6#7\relax{\def\x{#1#2#3#4#5#6}}%
\expandafter\x\fmtname xxxxxx\relax \def\y{splain}%
\ifx\x\y   % LaTeX or SliTeX?
\gdef\SetFigFont#1#2#3{%
  \ifnum #1<17\tiny\else \ifnum #1<20\small\else
  \ifnum #1<24\normalsize\else \ifnum #1<29\large\else
  \ifnum #1<34\Large\else \ifnum #1<41\LARGE\else
     \huge\fi\fi\fi\fi\fi\fi
  \csname #3\endcsname}%
\else
\gdef\SetFigFont#1#2#3{\begingroup
  \count@#1\relax \ifnum 25<\count@\count@25\fi
  \def\x{\endgroup\@setsize\SetFigFont{#2pt}}%
  \expandafter\x
    \csname \romannumeral\the\count@ pt\expandafter\endcsname
    \csname @\romannumeral\the\count@ pt\endcsname
  \csname #3\endcsname}%
\fi
\fi\endgroup
\begin{picture}(1329,1329)(439,-613)
\put(1756,-556){\makebox(0,0)[rb]{\smash{\SetFigFont{10}{12.0}{rm}{\color[rgb]{0,0,0}(1)}%
}}}
\put(1071, 29){\makebox(0,0)[rb]{\smash{\SetFigFont{8}{9.6}{rm}{\color[rgb]{0,0,0}2}%
}}}
\end{picture}
 \end {center}

(the balancing condition at e.g.\ the left vertex is $ 2 \cdot (1,1) + (-1,0) +
(-1,-2) = (0,0) $). The \emph {degree} of a tropical curve is the data of the
slopes of the unbounded edges together with their weights; we have to fix this
for enumerative questions. It is easily seen that the above tropical curve (1)
cannot be deformed (with its degree, i.e.\ the slopes of its unbounded edges
fixed) to another tropical curve that still passes through the three given
points. 

Let us now move one of the three given points, say the rightmost one down. This
has the effect of making the internal edge shorter as in (2) below until it
finally disappears and the graph acquires a vertex of valence 4 as in (3). But
if we move the point further down something strange happens: our tropical curve
now deforms into \emph {two different} curves (4) and (5) that pass through the
given points. So a naive count would give the result that the number of
tropical curves through the given points \emph {does depend} on the position of
the points, as we have deformed one curve (1) into two curves (4) and (5).

\begin {center} \begin{picture}(0,0)%
\includegraphics{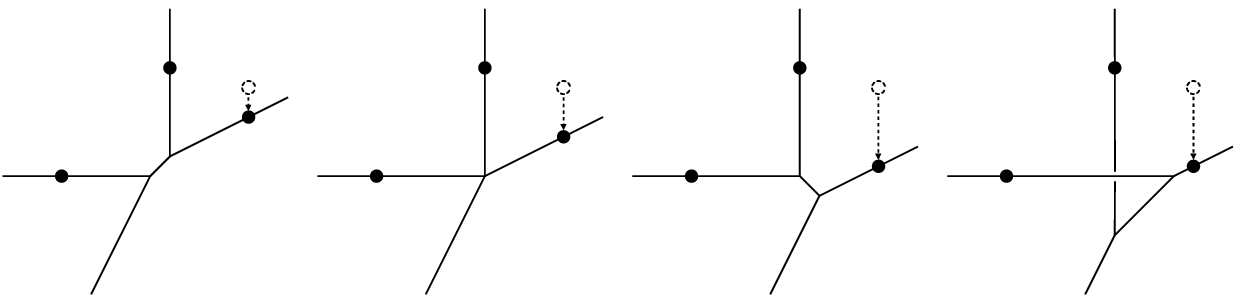}%
\end{picture}%
\setlength{\unitlength}{4144sp}%
\begingroup\makeatletter\ifx\SetFigFont\undefined
% extract first six characters in \fmtname
\def\x#1#2#3#4#5#6#7\relax{\def\x{#1#2#3#4#5#6}}%
\expandafter\x\fmtname xxxxxx\relax \def\y{splain}%
\ifx\x\y   % LaTeX or SliTeX?
\gdef\SetFigFont#1#2#3{%
  \ifnum #1<17\tiny\else \ifnum #1<20\small\else
  \ifnum #1<24\normalsize\else \ifnum #1<29\large\else
  \ifnum #1<34\Large\else \ifnum #1<41\LARGE\else
     \huge\fi\fi\fi\fi\fi\fi
  \csname #3\endcsname}%
\else
\gdef\SetFigFont#1#2#3{\begingroup
  \count@#1\relax \ifnum 25<\count@\count@25\fi
  \def\x{\endgroup\@setsize\SetFigFont{#2pt}}%
  \expandafter\x
    \csname \romannumeral\the\count@ pt\expandafter\endcsname
    \csname @\romannumeral\the\count@ pt\endcsname
  \csname #3\endcsname}%
\fi
\fi\endgroup
\begin{picture}(5649,1329)(439,-613)
\put(1756,-556){\makebox(0,0)[rb]{\smash{\SetFigFont{10}{12.0}{rm}{\color[rgb]{0,0,0}(2)}%
}}}
\put(3196,-556){\makebox(0,0)[rb]{\smash{\SetFigFont{10}{12.0}{rm}{\color[rgb]{0,0,0}(3)}%
}}}
\put(4636,-556){\makebox(0,0)[rb]{\smash{\SetFigFont{10}{12.0}{rm}{\color[rgb]{0,0,0}(4)}%
}}}
\put(6076,-556){\makebox(0,0)[rb]{\smash{\SetFigFont{10}{12.0}{rm}{\color[rgb]{0,0,0}(5)}%
}}}
\put(1156,-26){\makebox(0,0)[rb]{\smash{\SetFigFont{8}{9.6}{rm}{\color[rgb]{0,0,0}2}%
}}}
\end{picture}
 \end {center}

The solution to this problem is that tropical curves have to be counted with
multiplicities (as it is in fact required by the Correspondence Theorem). In
contrast to the complex case these multiplicities can be greater than 1 even
for points in general position. In the above example they turn out to be 4 for
(1), 3 for (4), and 1 for (5), so that the total weighted sum before and after
the deformation is still the same (for the precise definition of the
multiplicities see definition \ref{def-mult}).

In fact, this example shows already the general idea of our proof. In the space
$ \R^{2n} $ of $n$ points in the plane we consider the subspace of
configurations at which the topology of the tropical curves passing through
them changes, as e.g.\ in (3) above. This is a space of real codimension 1 in $
\R^{2n} $ (in fact it is a union of convex polyhedra of dimension $ 2n-1 $).
Whenever we pass such a ``wall'' in $ \R^{2n} $ we have to show that the
weighted sum of tropical curves through the points before and after crossing
the wall is the same. Once we have done this we know that we can go from any
configuration of points to any other along a path in $ \R^{2n} $, and that the
corresponding (weighted) count of tropical curves through these points will
remain constant along the path (in particular when crossing the walls). In the
whole analysis we can of course neglect the ``boundaries of the walls'' since
they have codimension 2 in $ \R^{2n} $ and can be avoided by the path. This
amounts to only consider configurations of points in sufficiently general
position.

This paper is organized as follows. In section \ref{sec-planetrop} we give a
rigorous definition of tropical curves with marked points. We then build up
suitable moduli spaces of such marked curves in section \ref{sec-moduli}.
Section \ref{sec-enum} studies the ``evaluation map'' that maps an $n$-marked
tropical curve to its configuration of points in $ \R^{2n} $. The main work of
our paper, i.e.\ the ``wall crossing analysis'' described above, is contained
in the proof of our main theorem \ref{constant}.

Finally we should say clearly that our paper gives a rather ad hoc solution to
the problem described above. From a theoretical point of view it would
certainly be more desirable to solve the problem in the same way as in complex
geometry: develop a tropical version of intersection theory and define the
numbers of curves through given points as intersection products on suitable
``moduli spaces of tropical stable maps''. Different collections of points
should then just correspond to ``rationally equivalent cycles'' in $ \R^{2n} $,
from which it should follow by general principles that the resulting
intersection products are the same. We hope that the ideas of our paper will be
useful to set up such a theory.

%%%%%%%%%%%%%%%%%%%%%%%%%%%%%%%%%%%%%%%%%%%%%%%%%%%%%%%%%%%%%%%%%%%%%%%%%%%%%%%

\section {Plane tropical curves} \label{sec-planetrop}

Let us start by defining tropical curves. Our definition differs slightly from
the ones given in \cite{Mi03} and \cite{NS04}. The reason for this is that
\cite{Mi03} and \cite{NS04} consider mostly ``generic tropical curves'', i.e.\
tropical curves that pass through given points (or subspaces) in ``general
position''. In our paper however the points will be in slightly less general
position, leading to ``degenerate'' tropical curves. As a consequence, we have
to be more careful in the definition of (possibly degenerate) tropical curves.

As in the classical case of complex geometry there are two ways to describe
curves in an ambient space: either as certain one-dimensional \emph {subspaces}
or as abstract curves together with a \emph {map} to the ambient space. While
these two notions describe the same objects for generic curves, they differ
e.g.\ in the cases where the curves have multiple components (resp.\ the map is
not generically injective). Both viewpoints have their advantages, and in fact
our definition of tropical curves in this paper is based on a mixture of these
two ideas. As a consequence the resulting moduli spaces will be quite well
suited for our computations, but unfortunately rather unsatisfactory from a
purely theoretic point of view.

Tropical curves are based on graphs. The graphs we will need are ``graphs with
multiple and unbounded edges allowed''. To make this precise and set up the
notation we will give a rigorous definition.

\begin {definition} \label{def-graph}
  A \df {graph} is a tuple $ \Gamma=(\Gamma^0,\Gamma',\partial,j) $ where
  \begin {itemize}
  \item $ \Gamma^0 $ and $ \Gamma' $ are finite sets (whose elements are
    called \df {vertices} and \df {flags}, respectively);
  \item $ \partial: \Gamma' \to \Gamma^0 $ is a map (called the \df {boundary
    map});
  \item $ j: \Gamma' \to \Gamma' $ is a map with $ j \circ j = \id
    $ (called the \df {glueing map}).
  \end {itemize}
  By the condition $ j \circ j = \id $ the relation
    \[ F_1 \sim F_2 :\Longleftrightarrow F_1 = F_2 \mbox { or } F_1 = j(F_2) \]
  for $ F_1,F_2 \in \Gamma' $ is an equivalence relation on the $ \Gamma' $
  whose equivalence classes all consist of one or two elements. The equivalence
  class $ [F] $ of a flag $ F \in \Gamma' $ is called an \df {edge} of
  $ \Gamma $. It is called an \df {unbounded edge} or \df {end} if $ [F] =
  \{F\} $ has only one element, and an \df {internal edge} otherwise. We denote
  the sets of edges, unbounded edges, and internal edges by $ \Gamma^1 $, $
  \Gamma^1_\infty $, and $ \Gamma^1_0 $, respectively. For a vertex $ V \in
  \Gamma^0 $ we denote by $ \val V := \#\{ F \in \Gamma' ;\; \partial F = V \}
  \in \N $ the \df {valence} of $V$.

  The \df {topological model} $ |\Gamma| $ of a graph $ \Gamma $ is obtained by
  starting with a disjoint union of a point $ \ast_V $ for every vertex $ V \in
  \Gamma^0 $ and a semi-open interval $ [0,1)_F $ for every flag $ F \in
  \Gamma' $, and glueing
  \begin {itemize}
  \item the point $ 0 \in [0,1)_F $ to the point $ \ast_{\partial F} $ for
    every flag $ F \in \Gamma' $;
  \item the intervals $ [0,1)_{F_1} $ and $ [0,1)_{F_2} $ on the common open
    subset $ (0,1) $ along the map $ t \mapsto 1-t $ for every pair of flags
    $ F_1,F_2 $ with $ F_1 \sim F_2 $ and $ F_1 \neq F_2 $.
  \end {itemize}
  A graph is called \df {connected} if its topological model is. The \df
  {genus} of a connected graph $ \Gamma $ is defined to be $ g(\Gamma) = \dim
  H^1 (|\Gamma|,\R) $.

  In the topological model every vertex $ V \in \Gamma^0 $ corresponds to a
  point $ \ast_V \in |\Gamma| $, every flag $ F \in \Gamma' $ corresponds to a
  semi-open interval $ [0,1) \subset |\Gamma| $, and every edge $ E \in
  \Gamma^1 $ (without its endpoints) corresponds to an open interval $ (0,1)
  \subset |\Gamma| $. By abuse of notation we will denote these points and
  intervals simply by $ V \in |\Gamma| $, $ F \subset |\Gamma| $, and
  $ E \subset |\Gamma| $, respectively.
\end {definition}

\begin {example} \label{ex-graph}
  The following picture shows (the topological model of) a graph $ \Gamma $ of
  genus 1 with 2 vertices $ \Gamma^0 = \{V_1,V_2\} $ (of valence 4 each) and 8
  flags $ \Gamma' = \{F_1,\dots,F_8\} $:

  \begin {center} \begin{picture}(0,0)%
\includegraphics{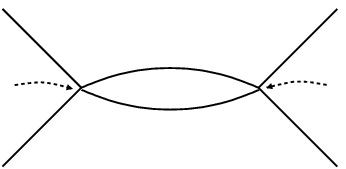}%
\end{picture}%
\setlength{\unitlength}{4144sp}%
\begingroup\makeatletter\ifx\SetFigFont\undefined
% extract first six characters in \fmtname
\def\x#1#2#3#4#5#6#7\relax{\def\x{#1#2#3#4#5#6}}%
\expandafter\x\fmtname xxxxxx\relax \def\y{splain}%
\ifx\x\y   % LaTeX or SliTeX?
\gdef\SetFigFont#1#2#3{%
  \ifnum #1<17\tiny\else \ifnum #1<20\small\else
  \ifnum #1<24\normalsize\else \ifnum #1<29\large\else
  \ifnum #1<34\Large\else \ifnum #1<41\LARGE\else
     \huge\fi\fi\fi\fi\fi\fi
  \csname #3\endcsname}%
\else
\gdef\SetFigFont#1#2#3{\begingroup
  \count@#1\relax \ifnum 25<\count@\count@25\fi
  \def\x{\endgroup\@setsize\SetFigFont{#2pt}}%
  \expandafter\x
    \csname \romannumeral\the\count@ pt\expandafter\endcsname
    \csname @\romannumeral\the\count@ pt\endcsname
  \csname #3\endcsname}%
\fi
\fi\endgroup
\begin{picture}(1554,744)(304,-433)
\put(541,174){\makebox(0,0)[b]{\smash{\SetFigFont{8}{9.6}{rm}{\color[rgb]{0,0,0}$F_1$}%
}}}
\put(516,-396){\makebox(0,0)[b]{\smash{\SetFigFont{8}{9.6}{rm}{\color[rgb]{0,0,0}$F_2$}%
}}}
\put(1636,-361){\makebox(0,0)[b]{\smash{\SetFigFont{8}{9.6}{rm}{\color[rgb]{0,0,0}$F_3$}%
}}}
\put(1671,189){\makebox(0,0)[b]{\smash{\SetFigFont{8}{9.6}{rm}{\color[rgb]{0,0,0}$F_4$}%
}}}
\put(841, 34){\makebox(0,0)[b]{\smash{\SetFigFont{8}{9.6}{rm}{\color[rgb]{0,0,0}$F_5$}%
}}}
\put(841,-256){\makebox(0,0)[b]{\smash{\SetFigFont{8}{9.6}{rm}{\color[rgb]{0,0,0}$F_6$}%
}}}
\put(1351,-256){\makebox(0,0)[b]{\smash{\SetFigFont{8}{9.6}{rm}{\color[rgb]{0,0,0}$F_7$}%
}}}
\put(1356, 29){\makebox(0,0)[b]{\smash{\SetFigFont{8}{9.6}{rm}{\color[rgb]{0,0,0}$F_8$}%
}}}
\put(334,-93){\makebox(0,0)[rb]{\smash{\SetFigFont{8}{9.6}{rm}{\color[rgb]{0,0,0}$V_1$}%
}}}
\put(1846,-97){\makebox(0,0)[lb]{\smash{\SetFigFont{8}{9.6}{rm}{\color[rgb]{0,0,0}$V_2$}%
}}}
\end{picture}
 \end {center}

  The graph has 6 edges: four unbounded edges $ E_i = \{F_i\} $ for $
  i=1,\dots,4 $, and two internal edges $ E_5 = \{ F_5,F_8 \} $, $ E_6 = \{
  F_6,F_7 \} $.
\end {example}

We are now ready to give the definition of a (plane) tropical curve.

\begin {definition} \label{def-tropcurve}
  A \df {tropical curve} is a triple $ C=(\Gamma,\omega,h) $, where
  \begin {enumerate}
  \item \label{def-tropcurve-a}
    $ \Gamma $ is a connected graph;
  \item \label{def-tropcurve-b}
    $ \omega: \Gamma^1 \to \N_{>0} $ is a map (called the \df {weight
    function});
  \item \label{def-tropcurve-c}
    $ h: |\Gamma| \to \R^2 $ is a continuous proper map that embeds every flag
    $ F \subset |\Gamma| $ into a (unique) affine line in $ \R^2 $ with
    rational slope. Moreover, if we denote by $ u(F) \in \Z^2 $ the primitive
    integral vector that starts at $ h(\partial F) $ and points in the
    direction of $ h(F) $ we set $ v(F) := \omega([F]) \cdot u(F) $ and
    require that at every vertex $ V \in \Gamma^0 $
    \begin {itemize}
    \item the vectors $ \{ v(F);\; F \in \Gamma' \mbox { with } \partial F=V \}
      $ span $ \R^2 $ as a vector space; and
    \item the \df {balancing condition}
        \[ \sum_{F \in \Gamma':\, \partial F = V}
             v(F) = 0 \]
      holds.
    \end {itemize}
  \end {enumerate}
  The genus $ g(C) $ of a tropical curve $ C=(\Gamma,\omega,h) $ is defined to
  be the genus of its underlying graph $ \Gamma $. Two tropical curves $
  (\Gamma,\omega,h) $ and $ (\tilde \Gamma, \tilde \omega, \tilde h) $ are
  called isomorphic if there is a homeomorphism $ \varphi: \Gamma \to \tilde
  \Gamma $ such that $ \tilde h \circ \varphi = h $ and $ \omega(E) = \tilde
  \omega(\varphi(E)) $ for all edges $ E \in \Gamma^1 $.
\end {definition}

\begin {remark} \label{rem-tropcurve}
  The condition in definition \ref{def-tropcurve} \ref{def-tropcurve-c} that
  $h$ is proper means precisely that the unbounded edges of $ \Gamma $ map to
  unbounded rays in $ \R^2 $. Note that the requirement that the flags around
  every vertex span $ \R^2 $ imply together with the balancing condition that
  every vertex has valence at least 3. If two flags around a vertex are mapped
  to the same line in $ \R^2 $ then this vertex must have valence at least 4.
\end {remark}

\begin {example} \label{ex-tropcurve}
  The following picture shows an example of a tropical curve $
  C=(\Gamma,\omega,h) $ of genus 1 based on the graph $ \Gamma $ of example
  \ref{ex-graph}. We have labeled the edges $ E \in \Gamma^1 $ by their weight
  $ \omega(E) $ and left out this label if the weight is 1.

  \begin {center} \begin{picture}(0,0)%
\includegraphics{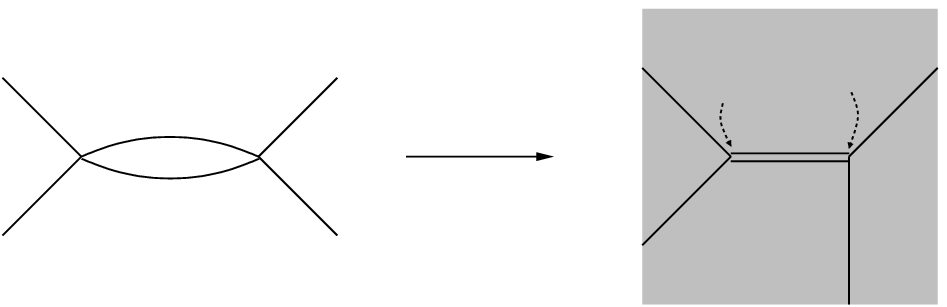}%
\end{picture}%
\setlength{\unitlength}{4144sp}%
\begingroup\makeatletter\ifx\SetFigFont\undefined
% extract first six characters in \fmtname
\def\x#1#2#3#4#5#6#7\relax{\def\x{#1#2#3#4#5#6}}%
\expandafter\x\fmtname xxxxxx\relax \def\y{splain}%
\ifx\x\y   % LaTeX or SliTeX?
\gdef\SetFigFont#1#2#3{%
  \ifnum #1<17\tiny\else \ifnum #1<20\small\else
  \ifnum #1<24\normalsize\else \ifnum #1<29\large\else
  \ifnum #1<34\Large\else \ifnum #1<41\LARGE\else
     \huge\fi\fi\fi\fi\fi\fi
  \csname #3\endcsname}%
\else
\gdef\SetFigFont#1#2#3{\begingroup
  \count@#1\relax \ifnum 25<\count@\count@25\fi
  \def\x{\endgroup\@setsize\SetFigFont{#2pt}}%
  \expandafter\x
    \csname \romannumeral\the\count@ pt\expandafter\endcsname
    \csname @\romannumeral\the\count@ pt\endcsname
  \csname #3\endcsname}%
\fi
\fi\endgroup
\begin{picture}(4422,1363)(304,-748)
\put(2476,-16){\makebox(0,0)[b]{\smash{\SetFigFont{8}{9.6}{rm}{\color[rgb]{0,0,0}$h$}%
}}}
\put(1451, 29){\makebox(0,0)[b]{\smash{\SetFigFont{8}{9.6}{rm}{\color[rgb]{0,0,0}$V_2$}%
}}}
\put(681, 29){\makebox(0,0)[b]{\smash{\SetFigFont{8}{9.6}{rm}{\color[rgb]{0,0,0}$V_1$}%
}}}
\put(4242,-484){\makebox(0,0)[b]{\smash{\SetFigFont{8}{9.6}{rm}{\color[rgb]{0,0,0}2}%
}}}
\put(4413, 29){\makebox(0,0)[b]{\smash{\SetFigFont{8}{9.6}{rm}{\color[rgb]{0,0,0}2}%
}}}
\put(1730,-234){\makebox(0,0)[b]{\smash{\SetFigFont{8}{9.6}{rm}{\color[rgb]{0,0,0}2}%
}}}
\put(1739, 45){\makebox(0,0)[b]{\smash{\SetFigFont{8}{9.6}{rm}{\color[rgb]{0,0,0}2}%
}}}
\put(3621,247){\makebox(0,0)[b]{\smash{\SetFigFont{8}{9.6}{rm}{\color[rgb]{0,0,0}$h(V_1)$}%
}}}
\put(4198,275){\makebox(0,0)[b]{\smash{\SetFigFont{8}{9.6}{rm}{\color[rgb]{0,0,0}$h(V_2)$}%
}}}
\put(1081,-556){\makebox(0,0)[b]{\smash{\SetFigFont{8}{9.6}{rm}{\color[rgb]{0,0,0}$\Gamma$}%
}}}
\put(4726,479){\makebox(0,0)[b]{\smash{\SetFigFont{8}{9.6}{rm}{\color[rgb]{0,0,0}$\R^2$}%
}}}
\end{picture}
 \end {center}

  In this example the two internal edges are mapped to the same image line.
  Note that the images of the flags around both vertices span $ \R^2 $ as a
  vector space as required by definition \ref{def-tropcurve}
  \ref{def-tropcurve-c}. Although $h$ embeds every flag in $ \R^2 $ it is not a
  global embedding. The balancing conditions at the two vertices are
  \begin {align*}
    (-1,1) + (-1,-1)+(1,0)+(1,0) &= (0,0) \mbox { at $ V_1 $} \\
    \mbox {and} \quad
    (-1,0)+(-1,0)+2\cdot(0,-1)+2\cdot(1,1) &=(0,0) \mbox { at $ V_2 $}.
  \end {align*}
  In the rest of this paper when we draw tropical curves we will for
  simplicity only draw the image $ h(\Gamma) $, with the edges labeled by their
  weight.
\end {example}

An important notion is that of the degree of a tropical curve:

\begin {definition} \label{def-degree}
  Let $G$ be the free abelian semigroup generated by $ \Z^2 \backslash
  \{(0,0)\} $. We denote the addition in $G$ by the symbol $ \oplus $ to
  distinguish it from the addition in $ \Z^2 $. For an element $ \Delta = u_1
  \oplus \cdots \oplus u_n $ that is a sum of $n$ vectors $ u_1,\dots,u_n \in
  \Z^2 $ we set $ \#\Delta := n $.

  If $ C=(\Gamma,\omega,h) $ is a tropical curve we define the \df {degree}
  of $C$ to be
    \[ \deg C := \bigoplus_{\{F\} \in \Gamma^1_\infty} v(F) \]
  where $ v(F) $ is as in definition \ref{def-tropcurve} \ref{def-tropcurve-c}.
\end {definition}

\begin {example} \label{ex-degree}
  The tropical curve $C$ of example \ref{ex-tropcurve} has degree
    \[ \deg C = (-1,1) \oplus (-1,-1) \oplus (0,-2) \oplus (2,2). \]
  Note that for every tropical curve $ C=(\Gamma,\omega,h) $ we have $ \#\deg C
  = \#\Gamma^1_\infty $ by definition. Moreover, if $ \deg C = v_1 \oplus
  \cdots \oplus v_n $ with $ v_1,\dots,v_n \in \Z^2 $ we must have $ v_1 +
  \cdots + v_n = 0 $ by the balancing condition.
\end {example}

To study tropical curves through given points in the plane we now have to
consider curves with marked points on them:

\begin {definition} \label{def-marked}
  For $ n \in \N $ we say that an \df {$n$-marked tropical curve} is a tuple $
  (C,x_1,\dots,x_n) $ where
  \begin {enumerate}
  \item \label{def-marked-a}
    $ C=(\Gamma,\omega,h) $ is a tropical curve;
  \item \label{def-marked-b}
    $ x_1,\dots,x_n \in |\Gamma| $ are points with the following property: for
    every $ i=1,\dots,n $ there is a flag $ F_i \subset |\Gamma| $ with $ x_i
    \in F_i $ such that $ |\Gamma| \backslash ([F_1] \cup \cdots \cup [F_n]) $
    contains no loops and no connected components with more than one unbounded
    end. (Here as usual $ [F] $ denotes the (open) edge corresponding to a flag
    $F$. Note that some of the points $ x_1,\dots,x_n $ may be vertices, and it
    is allowed that some of them coincide.)
  \end {enumerate}
  An isomorphism $ (C,x_1,\dots,x_n) \to (\tilde C,\tilde x_1,\dots,\tilde x_n)
  $ is an isomorphism $ C \to \tilde C $ of tropical curves taking $ x_i $ to $
  \tilde x_i $ for all $ i=1,\dots,n $.
\end {definition}

\begin {remark} \label{rem-marked}
  Let $ C=(\Gamma,\omega,h) $ be a tropical curve, and let $ x_1,\dots,x_n \in
  C $ be points none of which lies on a vertex of $ \Gamma $. Then in part
  \ref{def-marked-b} of definition \ref{def-marked} the edges $ [F_i] $ for the
  marked points $ x_i $ are uniquely defined: they are simply the edges on
  which the marked points lie. So in this case $ (C,x_1,\dots,x_n) $ is a
  marked tropical curve if and only if $ |\Gamma| \backslash \{x_1,\dots,x_n\}
  $ contains no loops and no connected components with more than one unbounded
  end.
\end {remark}

\begin {example} \label{ex-marked}
  In the picture below the left choice of 4 marked points makes the tropical
  curve of example \ref{ex-tropcurve} into a marked tropical curve, whereas the
  right one does not.

  \begin {center} \input {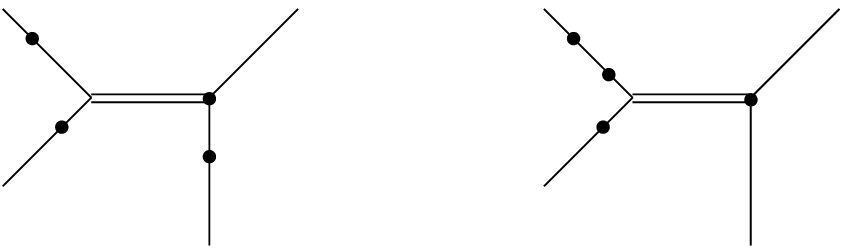} \end {center}

  To see that condition \ref{def-marked-b} of definition \ref{def-marked} is
  satisfied for the curve on the left but not for the one on the right note
  first that (analogously to remark \ref{rem-marked}) the edge $ [F_i] $
  corresponding to the point $ x_i $ is uniquely determined for $ i=1,2,3 $ in
  both cases: it is simply the edge on which $ x_i $ lies. For the flag $ F_4 $
  however we have four choices, namely the four flags starting at the vertex $
  x_4 $. For the curve on the left condition \ref{def-marked-b} of definition
  \ref{def-marked} is satisfied if we pick for $ F_4 $ one of the flags
  pointing to the left. For the curve on the right however any choice of flag $
  F_4 $ would lead to a connected component of the space $ |\Gamma| \backslash
  ([F_1] \cup \cdots \cup [F_4]) $ that has either two unbounded ends (if we
  pick a flag pointing to the left) or a loop (if we pick one of the other
  flags).
\end {example}

\begin {remark} \label{rem-marked-dim}
  It is easy to see that there is always a certain lower bound on the number
  $n$ of marked points that we need to make a given tropical curve of genus $g$
  and degree $ \Delta $ into a marked tropical curve. Recall that we have to
  make $ |\Gamma| $ into a space that has no loops and no connected component
  with more than one unbounded end by removing one (open) edge for every marked
  point. This means that we need at least $g$ marked points to break all loops
  in the graph and at least $ \#\Delta-1 $ more marked points to remove or
  separate all the unbounded ends. So we must always have the inequality $ n
  \ge \#\Delta+g-1 $ for an $n$-marked tropical curve of genus $g$ and degree $
  \Delta $.
\end {remark}

\begin {remark}
  The reason for condition \ref{def-marked-b} in definition \ref{def-marked} is
  that it ensures for sufficiently generic tropical curves that there are no
  deformations of the curve if we fix the images $ h(x_i) $ of the marked
  points in $ \R^2 $ (see proposition \ref{linear} for a precise statement).
  For example, in example \ref{ex-marked} the curve on the left cannot be
  deformed with fixed images of the marked points, whereas in the curve on the
  right the vertex on which $ x_4 $ lies can be moved to the right (thereby
  increasing the length of the ``double edge'', moving the two edges of weight
  2 horizontally, and letting $ x_4 $ move onto one of the middle edges).

  For more special tropical curves however it is unlikely that definition
  \ref{def-marked} does the ``right thing''. So if one wants to build up a
  general theory of marked tropical curves then definition \ref{def-marked}
  will probably not be the ``correct'' one.
\end {remark}

%%%%%%%%%%%%%%%%%%%%%%%%%%%%%%%%%%%%%%%%%%%%%%%%%%%%%%%%%%%%%%%%%%%%%%%%%%%%%%%

\section {Moduli spaces of tropical curves} \label{sec-moduli}

We are now ready to define the moduli spaces of marked tropical curves that
will be our main object of study. For the rest of this paper we fix the
following set-up: we choose a non-negative integer $g$ and a degree $ \Delta
\in G $, and we set $ n:= \#\Delta+g-1 $. This number $n$ (which has already
appeared in remark \ref{rem-marked-dim}) will turn out to be the required
number of points through which we get a finite non-zero number of tropical
curves of genus $g$ and degree $ \Delta $. We will therefore only consider
tropical curves with precisely this number $n$ of marked points.

\begin {definition} \label{def-moduli}
  The \df {moduli space $ \overline \calM_{g,\Delta} $ of marked tropical
  curves of genus $g$ and degree $ \Delta $} is defined to be the set of all
  $n$-marked tropical curves $ (C,x_1,\dots,x_n) $ of degree $ \Delta $ modulo
  isomorphisms such that $ n = \#\Delta+g-1 $, $ g(C)\le g $, and at least $
  g-g(C) $ of the marked points lie on vertices.
\end {definition}

\begin {remark}
  We have to allow curves of lower genus in our definition to ensure that the
  moduli spaces are ``closed under degenerations'' (see proposition
  \ref{boundary}). For an example consider again the left curve of example
  \ref{ex-marked}. If we shrink the length of the double edge to zero then
  the resulting curve will be of genus 0 (since we do not allow edges of length
  zero). Note also that by definition \ref{def-marked} \ref{def-marked-b} there
  must be a marked point on every loop of a tropical curve. So if we reduce the
  genus of a curve by a degeneration process that shrinks $k$ loops to a point
  as above then this will necessarily result in $k$ marked points lying on
  vertices after the degeneration. This is the reason for the condition on the
  marked points in definition \ref{def-moduli}.
\end {remark}

Let us now study the structure of these moduli spaces. The first thing we do is
to sort the elements of $ \overline \calM_{g,\Delta} $ according to their
``combinatorial type'':

\begin {definition} \label{def-combtype}
  Let $ (C,x_1,\dots,x_n) \in \overline \calM_{g,\Delta} $ be a marked
  tropical curve, where $ C=(\Gamma,\omega,h) $. For every $ i=1,\dots,n $ we
  denote by $ s(i) \in \Gamma^0 \amalg \Gamma^1 $ the unique stratum (vertex or
  edge) of $ \Gamma $ on which $ x_i $ lies. The \df {combinatorial type} of $
  (C,x_1,\dots,x_n) $ is the data $ \alpha=(\Gamma,\omega,u,s) $ (where $
  u: \Gamma' \to \Z^2 \backslash \{(0,0)\} $ was introduced in definition
  \ref{def-tropcurve} \ref{def-tropcurve-c}), i.e.\ it is given by the weighted
  graph, the direction of every edge of the graph in $ \R^2 $, and the
  information on which edges (or vertices) the marked points lie. The \df
  {codimension} of such a combinatorial type $ \alpha $ is defined to be
  \begin {align*}
    \codim \alpha &:= \sum_{V \in \Gamma^0} (\val V - 3) \\
         &\qquad + (g-g(C)) \\
         &\qquad + \#\{ i=1,\dots,n;\; s(i) \in \Gamma^0 \}.
  \end {align*}
  We denote by $ \calM_{g,\Delta}^\alpha $ the subset of $ \overline
  \calM_{g,\Delta} $ that corresponds to marked tropical curves of
  combinatorial type $ \alpha $. (We will see in proposition \ref
  {combtype-dim} and example \ref {ex-codim} that $ \codim \alpha $ is closely
  related, but not always equal to the codimension of $ \calM_{g,\Delta}^\alpha
  $ in $ \overline \calM_{g,\Delta}$.)
\end {definition}

\begin {remark} \label{rem-combtype}
  Note that the codimension of a combinatorial type as defined above is visibly
  the sum of three non-negative integers. In particular it is always a
  non-negative integer itself, and it is zero if and only if the curves of this
  type have only vertices of valence 3, have genus \emph {equal} to $g$, and
  have no marked points on vertices. If this is the case then we see moreover
  by the argument of remark \ref{rem-marked-dim} that $n$ is the \emph
  {minimal} number of marked points needed to satisfy the condition
  of remark \ref{rem-marked}. This means that then $ |\Gamma| \backslash
  \{x_1,\dots,x_n\} $ has exactly $ \#\Delta $ connected components, and that
  each of these components contains exactly one unbounded end.
\end {remark}

\begin {remark}\label{combtype-formula}
  For future computations we will need another description of the codimension
  of a combinatorial type $ \alpha $. To derive it, note that the genus of a
  graph $ \Gamma $ is given by $ g(C) = \#\Gamma^1_0 - \#\Gamma^0 + 1 $.
  Moreover, counting the vertices of the graph leads to the equation $
  \#\Gamma^1_\infty + 2 \#\Gamma^1_0 = \sum_{V \in \Gamma^0} \val V $.
  Combining these two equations and using the relations $ \#\Gamma^1_\infty =
  \#\Delta $ and $ n = \#\Delta + g - 1 $ we arrive at the result
    \[ \codim \alpha = 2n-2+2g(C)-\#\Gamma^1_0-\#\{i;\; s(i) \in \Gamma^1 \}.
       \]
\end {remark}

\begin {remark} \label{genus-2}
  Let $ \alpha $ be a combinatorial type occurring in $ \overline
  \calM_{g,\Delta} $ that corresponds to tropical curves of genus strictly less
  than $g$. Then there must be at least one marked point on a vertex by
  definition \ref{def-moduli}. In particular we must then have $ \codim \alpha
  \ge 2 $ by definition \ref{def-combtype}.
\end {remark}

\begin {proposition} \label{combtype-finite}
  There are only finitely many combinatorial types occurring in a given moduli
  space $ \overline \calM_{g,\Delta} $.
\end {proposition}

\begin {proof}
  This follows essentially from \cite{NS04} proposition 2.1. There it is shown
  that the number of possibilities for $(\Gamma, \omega, u)$ is finite for a
  given degree $ \Delta $ and genus of the curves. Furthermore, the genus of
  the curves in $ \overline \calM_{g,\Delta} $ is bounded by $g$, and there are
  only finitely many possibilities for $s$. 
\end {proof}

We will now study the spaces $ \calM_{g,\Delta}^\alpha $ for fixed
combinatorial types separately. Of course the idea of the codimension of a
combinatorial type $ \alpha $ is that it should correspond to the codimension
of its moduli space $ \calM_{g,\Delta}^\alpha $ in the whole space $ \overline
\calM_{g,\Delta} $. Unfortunately this is not true in general (see example
\ref{ex-codim}). It is true however in codimensions up to 2 (which will be
sufficient for our purposes) \emph {except} for one special case:

\begin {definition}
  We say that a combinatorial type $ \alpha=(\Gamma,\omega,u,s) $ is \df
  {exceptional} if $ \codim \alpha=2 $ and $ \Gamma $ has two vertices of
  valence 4 that are joined by two edges. In other words we can say that $
  \alpha $ is exceptional if and only if it contains the picture of example
  \ref{ex-graph} as a ``subgraph'', all other vertices are of valence 3, the
  genus of the graph is maximal, and no marked points lie on vertices.
\end {definition}

\begin {proposition} \label{combtype-dim}
  For every combinatorial type $ \alpha $ occurring in $ \overline
  \calM_{g,\Delta} $ the space $ \calM_{g,\Delta}^\alpha $ is naturally an
  (unbounded) open convex polyhedron in a real affine space, i.e.\ a subset of
  a real affine space given by finitely many linear strict inequalities. For
  its dimension we have
    \[ \dim \calM_{g,\Delta}^\alpha \;\; \begin {cases}
         \;\; = 2n & \mbox {if $ \codim \alpha = 0 $}; \\
         \;\; = 2n-1 & \mbox {if $ \codim \alpha = 1 $ or $ \alpha $ is
                  exceptional}; \\
         \;\; \le 2n-2 & \mbox {otherwise}.
       \end {cases} \]
\end {proposition}

\begin {proof}
  The proof of this proposition is based on the ideas of \cite {Mi03}
  proposition 2.23. Our result is similar to \cite {Shu04} lemma 2.2 but
  differs in that we consider parametrized and not embedded tropical curves (so
  that we cannot apply Shustin's technique of Newton polyhedra).

  Fixing the combinatorial type of a tropical curve $ C=(\Gamma,\omega,h) $
  simply means that we fix the weighted graph, the slopes of the images of all
  edges in $ \R^2 $, and the edges (resp.\ vertices) on which the marked points
  lie. In contrast, the combinatorial type does not fix the position of the
  curve in the plane, the lengths of the images of the internal edges, and the
  position of the marked points (that do not lie on vertices) on their
  respective edges. So $ \calM_{g,\Delta}^\alpha $ can be thought of as a
  subset of the real affine space $A$ whose coordinates are
  \begin {enumerate}
  \item \label{combtype-dim-a}
    the position $ h(V) \in \R^2 $ of a fixed ``root vertex'' $ V \in
    \Gamma^0 $;
  \item \label{combtype-dim-b}
    the lengths of the images $ h(E) \subset \R^2 $ of the internal edges $
    E \in \Gamma^1_0 $;
  \item \label{combtype-dim-c}
    for every marked point $ x_i $ lying on an edge $ s(i)=E_i \in \Gamma^1 $
    its distance in $ \R^2 $ from a neighboring vertex, i.e.\ the number $
    |h(x_i)-h(\partial F_i)| $ for a fixed flag $ F_i \in \Gamma' $ with $ F_i
    \in E_i $.
  \end {enumerate}
  Note that a different choice of root vertex in \ref{combtype-dim-a} or flag
  in \ref{combtype-dim-c} would simply correspond to an affine isomorphism.
  Therefore $ \calM_{g,\Delta}^\alpha $ is naturally a subset of a real affine
  space $A$. Moreover, by remark \ref{combtype-formula} the dimension of $A$ is
    \[ \dim A = 2+\#\Gamma^1_0+\#\{i;\; s(i) \in \Gamma^1 \}
         = 2n+2g(C)-\codim \alpha. \]
  Note however that the affine coordinates described above can only be chosen
  independently if $ g(C)=0 $. Otherwise, each loop in the graph $ \Gamma $
  leads to two linear equations on the lengths of its internal edges describing
  the condition that the image of this loop closes up in the plane $ \R^2 $. As
  it suffices to consider these conditions for a chosen set of generators of $
  H^1(\Gamma,\Z) $ we arrive at a total of $ 2g(C) $ linear conditions. So we
  would expect $ \calM_{g,\Delta}^\alpha $ to have dimension $ \dim A - 2g(C) =
  2n-\codim \alpha $. This is in general only a lower bound however since the $
  2g(C) $ conditions above need not be independent (see e.g.\ example
  \ref{ex-codim}). The main work in the proof of our proposition is now to give
  a good estimate on how many of these conditions are independent.

  To do so, pick a fixed tropical curve in $ \calM_{g,\Delta}^\alpha $ and a
  vertex $ V_1 \in \Gamma^0 $ of maximal valence. Let $ L \subset \R^2 $ be a
  fixed line through its image point $ h(V_1) $. (In some cases we will specify
  later which line is to be picked here. For the moment all our constructions
  work for an arbitrary line $L$.)

  Order the vertices of $ \Gamma $ starting with $ V_1 $ so that the distance
  of their image points from $L$ is increasing (if for some vertices this
  distance is equal we order them arbitrarily). Orient the edges so that they
  point from the lower to the higher vertex. The unbounded edges are always
  oriented so that they point from its vertex to infinity. Note that the
  balancing condition of definition \ref{def-tropcurve} \ref{def-tropcurve-c}
  implies that every vertex has at least one adjacent edge pointing in a
  direction of strictly increasing distance from $L$. So every vertex must have
  at least one adjacent edge which is oriented away from it. (This is in fact
  the only property of the chosen order and orientations that we will need.)
  An example is shown in the picture below on the left:

  \begin {center} \input {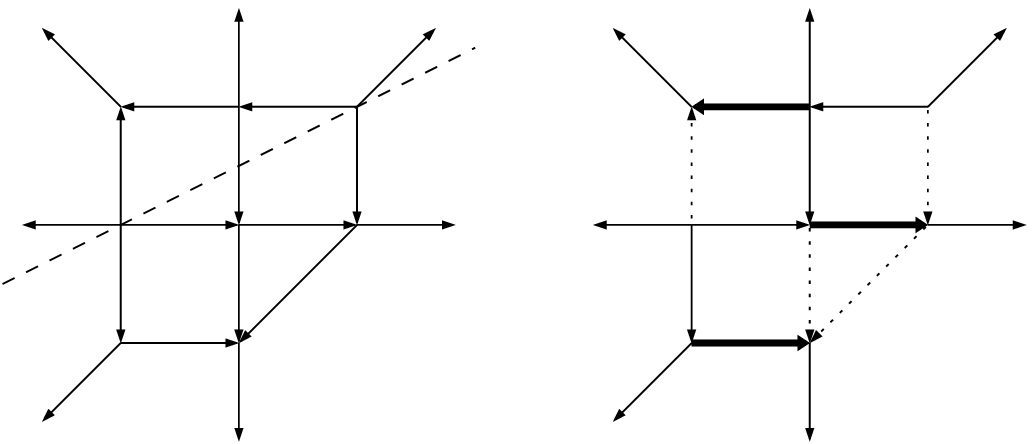} \end {center}

  We will now distinguish recursively $ 2g(C) $ edges $ E_1,\dots,E_{g(C)},
  E'_1,\dots,E'_{g(C)} $ as follows. For $ i=1,\dots,g(C) $ we let $ E_i $ be
  an (internal) edge contained in a loop of $ \Gamma \backslash \{
  E_1,\dots,E_{i-1} \} $ such that the vertex that this edge points to is
  the highest possible (in the chosen ordering). Then $ T := \Gamma \backslash
  \{ E_1,\dots,E_{g(C)} \} $ is a maximal tree in $ \Gamma $. In particular,
  for all $ i=1,\dots,g(C) $ the edge $ E_i $ closes a unique loop $ \Gamma_i $
  in $ T \cup \{E_i\} \subset \Gamma $. We let $ E'_i $ be the unique
  (internal) edge of $T$ that is contained in $ \Gamma_i $ and adjacent to the
  vertex that $ E_i $ points to. As an example, in the picture above on the
  right (where $ g(C)=4 $) we have drawn the edges $ E_1,\dots,E_4 $ as dotted
  arrows and $ E'_1,\dots,E'_4 $ in bold. The maximal tree $T$ consists exactly
  of all solid lines. Note that by construction the edges $ E_1,\dots,E_{g(C)}
  $ are all distinct and different from the $ E'_1,\dots,E'_{g(C)} $. It may
  happen however that not all edges $ E'_1,\dots,E'_{g(C)} $ are distinct. Note
  that by construction $ E_i $ and $ E'_i $ always point to the same vertex,
  namely to the highest vertex contained in the loop $ \Gamma_i $.

  We will now define a set of \emph{conditional edges} by starting with the
  $ 2g(C) $ edges $ E_1,\dots,E_{g(C)},E'_1,\dots,E'_{g(C)} $ and removing some
  of these edges by applying the following rules at each vertex $V$:
  \begin {itemize}
  \item[(i)] if there is at least one edge $ E'_i $ pointing to $V$ that is not
    parallel to its corresponding edge $ E_i $ then we keep the edge $ E'_i $
    with this property such that $i$ is maximal and remove all other edges $
    E'_1,\dots,E'_{g(C)} $ that point to $V$;
  \item [(ii)] if there is no such edge then we remove all edges $
    E'_1,\dots,E'_{g(C)} $ that point to $V$.
  \end {itemize}
  Note that all edges $ E_1,\dots,E_{g(C)} $ will end up to be conditional
  edges, and that all conditional edges will be distinct. In the example above
  we end up with the 7 conditional edges $ E_1,E_2,E_3,E_4,E'_2,E'_3,E'_4 $.

  We claim that for any (marked) tropical curve in $ \calM_{g,\Delta}^\alpha $
  the lengths of its conditional edges are determined uniquely in terms of the
  lengths of all other edges. To see this apply the following procedure
  recursively for $ i=g(C),\dots,1 $: assume that we know already the lengths
  of all edges $ E_{i+1},\dots,E_{g(C)},E'_{i+1},\dots,E'_{g(C)} $ as well as
  of all unconditional edges. Then by construction the only edges in the loop
  $ \Gamma_i $ whose lengths are not yet known can be $ E_i $ and $ E'_i $
  (if $V$ is the vertex that $ E_i $ and $ E'_i $ point to then all other edges
  in $ \Gamma_i $ must point to smaller vertices than $V$ whereas all edges $
  E_j $ and $ E'_j $ with $ j<i $ point to vertices greater than or equal to
  $V$). If $ E_i $ and $ E'_i $ are not parallel then the condition that $
  \Gamma_i $ closes up in $ \R^2 $ determines both their lengths uniquely.
  Otherwise $ E'_i $ is an unconditional edge by (ii), and $ E_i $ is again
  determined uniquely by the condition that $ \Gamma_i $ closes up.

  It follows that the dimension of $ \calM_{g,\Delta}^\alpha $ is at most equal
  to $ \dim A $ minus the number of conditional edges. So let us determine how
  many conditional edges there are. Note that when going from the $ 2g(C) $
  edges $ E_1,\dots,E_{g(C)},E'_1, \dots,E'_{g(C)} $ to the conditional edges
  we removed at most $ \val V - 3 $ edges at each vertex $V$: in case (i) above
  there is at least one edge pointing away from $V$ and one pair $ \{ E_i,E'_i
  \} $ that we do not remove. In case (ii), if we remove any edge $ E'_i $ at
  all there is at least one edge pointing away from $V$ and one other edge $
  E_i $ that we do not remove. But these cannot be all edges adjacent to $V$
  since then the flags adjacent to $V$ would not span $ \R^2 $ in contrast to
  definition \ref{def-tropcurve} \ref{def-tropcurve-c}. So there must be at
  least one other edge that is not removed, leading again to a total of at
  least three edges at $V$ that are not removed.

  Keeping in mind that we do not remove any edge at the vertex $ V_1 $ at all
  (since no edge points towards it) it follows that the number of conditional
  edges is at least
  \begin {align*}
    & 2g(C)-\sum_{V \neq V_1} (\val V-3) \\
    & \qquad = 2g(C)-\codim \alpha+(\val V_1-3)
      +(g-g(C))+\#\{i;\; s(i) \in \Gamma^0 \}
  \end {align*}
  so that an upper bound for the dimension of $ \calM_{g,\Delta}^\alpha $ is
    \[ 2n-(\val V_1-3)-(g-g(C))-\#\{i;\; s(i) \in \Gamma^0 \}. \tag {$*$} \]
  We now consider several cases, stopping at the first one that applies to
  $ \alpha $:
  \begin {itemize}
  \item If $ \codim \alpha\le 1 $ then $ V_1 $ is the only vertex that can
    possibly have valence greater than 3. So the number $ (*) $ is simply $
    2n-\codim \alpha $, and it follows that $ \dim \calM_{g,\Delta}^\alpha =
    2n-\codim \alpha $ (as we know already that this number is also a lower
    bound).
  \item If the number $(*)$ is at most $ 2n-2 $ then the dimension statement
    of the proposition follows immediately.
  \item If there are two vertices of valence 4 that are not joined by more than
    one edge then we choose $L$ above to be a line through these two vertices.
    We can then label these two vertices as $ V_1 $ and $ V_2 $. It follows
    that there is at most one edge pointing to $ V_2 $, i.e.\ that there are no
    edges removed at $ V_2 $ either in the above procedure. Hence we can
    subtract $ \val V_2-3 $ from the number $(*)$. So again we conclude that $
    \dim \calM_{g,\Delta}^\alpha \le 2n-2 $.
  \item The only case left is that we have $ g=g(C) $, no marked points on
    vertices, only vertices of valence 3 and 4, and all vertices of valence 4
    joined by at least 2 edges between each pair of them. This is only possible
    if there are exactly 2 vertices of valence 4, and if these vertices are
    joined by exactly 2 edges. In other words, $ \alpha $ is exceptional. In
    this case the dimension is obviously the same as for the combinatorial type
    where the double edge is replaced by one edge of added weight. As this new
    type has codimension 1 it follows that $ \dim \calM_{g,\Delta}^\alpha =2n-1
    $ in this case.
  \end {itemize}
  We have therefore proven the dimension statement of the proposition. So
  to finish the proof it suffices to notice that in the linear subspace of $A$
  determined by the $ 2g(C) $ (not necessarily independent) loop conditions the
  subset $ \calM_{g,\Delta}^\alpha $ is simply the open convex polyhedron
  given by the strict inequalitites that all lengths of the internal edges
  \ref{combtype-dim-b} are positive, and that all marked points on edges
  \ref{combtype-dim-c} lie in the interior of their respective edges. The
  polyhedron is unbounded since the coordinates of the root vertex
  \ref{combtype-dim-a} are not restricted.
\end {proof}

\begin {example} \label{ex-codim}
  The following example shows that we do \emph {not} have $ \dim
  \calM_{g,\Delta}^\alpha = 2n-\codim \alpha $ in general. We consider marked
  tropical curves of genus $ g=3 $ and degree $ \Delta = (-4,-2) \oplus (4,-2)
  \oplus (0,4) $ that have $ n = \#\Delta+g-1 = 5 $ marked points and are of
  the following combinatorial type $ \alpha $ (with an arbitrary choice of
  labeling of the marked points):

  \begin {center} \begin{picture}(0,0)%
\includegraphics{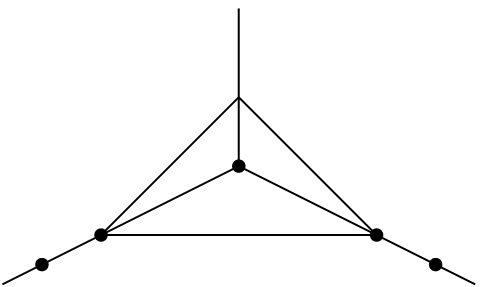}%
\end{picture}%
\setlength{\unitlength}{4144sp}%
\begingroup\makeatletter\ifx\SetFigFont\undefined
% extract first six characters in \fmtname
\def\x#1#2#3#4#5#6#7\relax{\def\x{#1#2#3#4#5#6}}%
\expandafter\x\fmtname xxxxxx\relax \def\y{splain}%
\ifx\x\y   % LaTeX or SliTeX?
\gdef\SetFigFont#1#2#3{%
  \ifnum #1<17\tiny\else \ifnum #1<20\small\else
  \ifnum #1<24\normalsize\else \ifnum #1<29\large\else
  \ifnum #1<34\Large\else \ifnum #1<41\LARGE\else
     \huge\fi\fi\fi\fi\fi\fi
  \csname #3\endcsname}%
\else
\gdef\SetFigFont#1#2#3{\begingroup
  \count@#1\relax \ifnum 25<\count@\count@25\fi
  \def\x{\endgroup\@setsize\SetFigFont{#2pt}}%
  \expandafter\x
    \csname \romannumeral\the\count@ pt\expandafter\endcsname
    \csname @\romannumeral\the\count@ pt\endcsname
  \csname #3\endcsname}%
\fi
\fi\endgroup
\begin{picture}(2184,1284)(439,-748)
\put(766,-736){\makebox(0,0)[b]{\smash{\SetFigFont{8}{9.6}{rm}{\color[rgb]{0,0,0}2}%
}}}
\put(2296,-736){\makebox(0,0)[b]{\smash{\SetFigFont{8}{9.6}{rm}{\color[rgb]{0,0,0}2}%
}}}
\put(1596,-146){\makebox(0,0)[b]{\smash{\SetFigFont{8}{9.6}{rm}{\color[rgb]{0,0,0}2}%
}}}
\put(1596,344){\makebox(0,0)[b]{\smash{\SetFigFont{8}{9.6}{rm}{\color[rgb]{0,0,0}4}%
}}}
\end{picture}
 \end {center}

  We have $ \codim \alpha = 6 $ since there are three 4-valent vertices in the
  graph and three marked points on vertices. We would therefore expect $
  \dim \calM_{g,\Delta}^\alpha $ to be $ 2n-\codim \alpha = 10-6 = 4 $. We have
  $ \dim \calM_{g,\Delta}^\alpha = 5 $ however: there are 2 dimensions for
  moving the marked points on edges, 2 dimensions for translations of the curve
  in the plane, and 1 more dimension for rescaling the whole curve.
\end {example}

\begin {definition}
  Let $ \alpha $ be a combinatorial type occurring in a moduli space $
  \overline \calM_{g,\Delta} $. By proposition \ref{combtype-dim} the space $
  \calM_{g,\Delta}^\alpha $ of curves of this given type is naturally an open
  subset of a real affine space $ A_\alpha $. We denote by $ \overline
  \calM_{g,\Delta}^\alpha $ the closure of $ \calM_{g,\Delta}^\alpha $ in $
  A_\alpha $.
\end {definition}

\begin {proposition} \label{boundary}
  Let $ \alpha $ be a combinatorial type occurring in a moduli space $
  \overline \calM_{g,\Delta} $. Then every point in $ \overline \calM_{g,
  \Delta}^\alpha $ can naturally be thought of as a marked tropical curve in $
  \overline \calM_{g,\Delta} $. The corresponding map $ i_\alpha: \overline
  \calM_{g,\Delta}^\alpha \to \overline \calM_{g,\Delta} $ maps the boundary $
  \partial \calM_{g,\Delta}^\alpha $ to the union of the strata $ \calM_{g,
  \Delta}^{\alpha'} $ such that
  \begin {enumerate}
  \item \label{boundary-a}
    the number of internal edges plus the number of marked points lying on
    edges is smaller for $ \alpha' $ than for $ \alpha $;
  \item \label{boundary-b}
    $ \codim \alpha' > \min (1,\codim \alpha) $.
  \end {enumerate}
  Moreover, the restriction of $ i_\alpha $ to any inverse image of such a
  stratum $ \calM_{g,\Delta}^{\alpha'} $ is an affine map.
\end {proposition}

\begin {proof}
  Recall that by the proof of proposition \ref{combtype-dim} the boundary
  $ \partial \calM_{g,\Delta}^{\alpha} $ is given by tuples $ (C,x_1,\dots,x_n)
  $ with $ C=(\Gamma,\omega,h) $ such that some marked points $ x_i $ with $
  s(i) \in \Gamma^1 $ lie on the boundary of their respective edges (i.e.\ on a
  vertex) and/or some interior edges are mapped to a line of length 0 in $ \R^2
  $ (i.e.\ to a point). In the former case this simply changes the
  combinatorial type $ (C,x_1,\dots,x_n) $ so that the points $ x_i $ are now
  required to map to a vertex instead of to an edge, i.e.\ so that $ s(i) \in
  \Gamma^0 $. In the latter case we will simply change the graph by removing
  all edges $E$ such that $ h(E) $ is a point and glueing the vertices in
  $\partial E$ to one vertex as in the following picture:

  \begin {center} \begin{picture}(0,0)%
\includegraphics{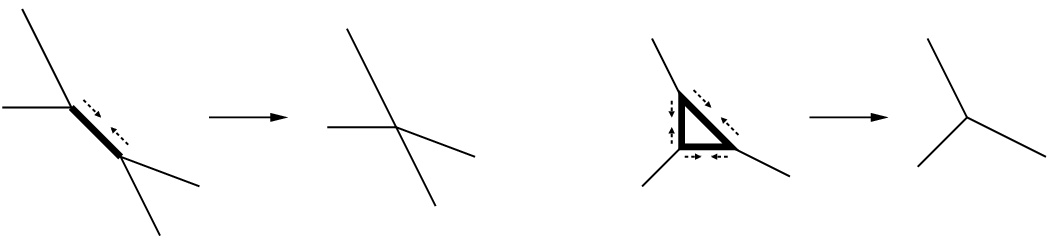}%
\end{picture}%
\setlength{\unitlength}{4144sp}%
\begingroup\makeatletter\ifx\SetFigFont\undefined
% extract first six characters in \fmtname
\def\x#1#2#3#4#5#6#7\relax{\def\x{#1#2#3#4#5#6}}%
\expandafter\x\fmtname xxxxxx\relax \def\y{splain}%
\ifx\x\y   % LaTeX or SliTeX?
\gdef\SetFigFont#1#2#3{%
  \ifnum #1<17\tiny\else \ifnum #1<20\small\else
  \ifnum #1<24\normalsize\else \ifnum #1<29\large\else
  \ifnum #1<34\Large\else \ifnum #1<41\LARGE\else
     \huge\fi\fi\fi\fi\fi\fi
  \csname #3\endcsname}%
\else
\gdef\SetFigFont#1#2#3{\begingroup
  \count@#1\relax \ifnum 25<\count@\count@25\fi
  \def\x{\endgroup\@setsize\SetFigFont{#2pt}}%
  \expandafter\x
    \csname \romannumeral\the\count@ pt\expandafter\endcsname
    \csname @\romannumeral\the\count@ pt\endcsname
  \csname #3\endcsname}%
\fi
\fi\endgroup
\begin{picture}(4794,1059)(529,-973)
\put(926,-606){\makebox(0,0)[b]{\smash{\SetFigFont{8}{9.6}{rm}{\color[rgb]{0,0,0}$2$}%
}}}
\end{picture}
 \end {center}

  In the picture we have drawn in bold the edges whose lengths tend to zero and
  are finally mapped to a point and removed.

  It is easy to see that the conditions of definition \ref{def-tropcurve}
  \ref{def-tropcurve-c} are still satisfied for the new graph. Hence the result
  is a tropical curve. Together with the marked points it is in fact an
  $n$-marked tropical curve (i.e.\ it satisfies condition \ref{def-marked-b}
  of definition \ref{def-marked}) since with the notation of definition
  \ref{def-marked} \ref{def-marked-b} none of the processes above can open
  up a loop in $ |\Gamma| \backslash ([F_1] \cup \cdots \cup [F_n]) $ or
  connect two unbounded edges in $ |\Gamma| \backslash ([F_1] \cup \cdots \cup
  [F_n]) $ that have not been connected before. This shows that the points in
  the boundary $ \partial \calM_{g,\Delta}^\alpha $ can naturally be thought of
  as marked tropical curves in $ \overline \calM_{g,\Delta} $ themselves.

  Let us now analyze which combinatorial types can occur on the boundary, i.e.\
  in the image $ i_\alpha(\partial \calM_{g,\Delta}^\alpha) $. It is clear by
  the construction above that condition \ref{boundary-a} of the proposition
  must be satisfied. To show \ref{boundary-b} we distinguish two cases: if a
  combinatorial type $ \alpha' $ occurring in the boundary has genus strictly
  less than $g$ then $ \codim \alpha' > 1 $ by remark \ref{genus-2}. On the
  other hand, if $ \alpha' $ (and hence also $ \alpha $) corresponds to curves
  of genus $g$ then we have $ \codim \alpha' > \codim \alpha $ by
  \ref{boundary-a} and the formula of remark \ref{combtype-formula}. This
  proves \ref{boundary-b}.
  
  Finally, it is clear that the restriction of $ i_\alpha $ to the inverse
  image of any stratum $ \calM_{g,\Delta}^{\alpha'} $ is an affine map since
  the affine structure on any stratum is given by the position of the curve in
  the plane, the lengths of the internal edges, and the position of the marked
  points on edges.
\end {proof}

\begin {remark} \label{rem-codim}
  One would expect that in general the boundary of a stratum $ \calM_{g,\Delta}
  ^\alpha $ corresponds only to marked curves of higher codimensions. This is
  not true however: the curve in example \ref{ex-codim} can be degenerated by
  shrinking the triangles to zero size, arriving at a tropical curve of genus 0
  with three edges and one vertex, and three of the marked points lying on the
  vertex. This new combinatorial type in the boundary has codimension 6,
  which is the same codimension as the one of the curves that we started with.
  Proposition \ref{boundary} \ref{boundary-b} ensures however that the
  codimension of the curves in the boundary is always bigger if we start with a
  stratum of codimension 0 or 1.
\end {remark}

\begin {remark} \label{moduli-topological}
  By propositions \ref{combtype-finite} and \ref{boundary} we can think of the
  moduli space $ \overline \calM_{g,\Delta} $ as being obtained by starting
  with finitely many unbounded closed convex polyhedra $ \overline
  \calM_{g,\Delta}^\alpha $ and then glueing them together by attaching each
  boundary $ \partial \calM_{g,\Delta}^\alpha $ with affine maps to some
  polyhedra $ \overline \calM_{g,\Delta}^{\alpha'} $ such that the number of
  internal edges plus the number of marked points lying on edges is smaller for
  $ \alpha' $ than for $ \alpha $. In particular, this makes $ \overline
  \calM_{g,\Delta} $ into a topological space with a natural stratification
  such that each stratum is an unbounded open convex polyhedron.

  One would expect that the moduli space $ \overline \calM_{g,\Delta} $ (or
  more probably a similar space with a slightly different definition of
  marked tropical curves) can in fact be given the structure of a tropical
  variety itself. However, the theory of abstract tropical varieties is still
  very much in its beginnings (see e.g.\ \cite {Mi05}) so that we will not use
  this language here.
\end {remark}

%%%%%%%%%%%%%%%%%%%%%%%%%%%%%%%%%%%%%%%%%%%%%%%%%%%%%%%%%%%%%%%%%%%%%%%%%%%%%%%

\section {Enumerative geometry of tropical curves} \label{sec-enum}

We are now ready to count tropical curves through some given points in the
plane. As in the previous section we fix a genus $ g \in \N $ and degree $
\Delta \in G $ and consider tropical curves with $ n := \#\Delta+g-1 $ marked
points.

\begin {definition} \label{def-eval}
  We define the \df {evaluation map} to be
  \begin {align*}
    \pi: \qquad \qquad \qquad \quad \overline \calM_{g,\Delta} &\to \R^{2n} \\
      ((\Gamma,\omega,h),x_1,\dots,x_n) &\mapsto (h(x_1),\dots,h(x_n)).
  \end {align*}
\end {definition}

\begin {proposition} \label{linear}
  The evaluation map $ \pi: \overline \calM_{g,\Delta} \to \R^{2n} $ is
  affine on each closed stratum $ \overline \calM_{g,\Delta}^\alpha $ (and
  hence in particular continuous). Moreover, it is injective on each stratum
  $ \calM_{g,\Delta}^\alpha $ with $ \codim \alpha = 0 $.
\end {proposition}

\begin {proof}
  For a fixed combinatorial type it is clear that the positions $ h(x_i) \in
  \R^2 $ of all marked points are linear functions in the affine coordinates
  constructed in the proof of proposition \ref{combtype-dim}. Hence $ \pi $ is
  affine on each closed stratum $ \overline \calM_{g,\Delta}^\alpha $.

  Now let $ \alpha $ be a combinatorial type with $ \codim \alpha=0 $, and let
  $ p_1,\ldots,p_n \in \R^2 $. We are going to show that $ \pi $ is injective
  on $ \calM_{g,\Delta}^\alpha $, i.e.\ that there is at most one marked
  tropical curve $ (C=(\Gamma,\omega,h),x_1,\ldots,x_n) $ of type $ \alpha $
  with $ h(x_i)=p_i $ for all $i$. Note that $ \codim \alpha = 0 $ implies in
  particular that all vertices of $ \Gamma $ have valence 3.

  As $\codim \alpha=0$ none of the marked points can lie on a vertex. Let $
  E_i =s(i) \in \Gamma^1 $ be the edge on which $ x_i $ lies. By remark
  \ref{rem-combtype} $ |\Gamma| \backslash \{x_1,\dots,x_n\} $ has $ \#\Delta $
  connected components each of which has genus 0 and contains exactly one
  unbounded end. Let $K$ be one of these connected components. We can pick two
  distinct points $ x_i $ and $ x_j $ in $ \overline K $ such that the
  corresponding edges $ E_i $ and $ E_j $ are adjacent to the same vertex $V$
  of $K$. Then both the direction and a starting point of the line segments $
  h(E_i) $ and $ h(E_j) $ are fixed by $u$ and the points $ p_i $ and $ p_j $.
  Note that the directions of the two line segments cannot be the same by
  remark \ref{rem-tropcurve} since $ \val V = 3 $. Hence these data determine
  $ h(V) $ (and thus the lengths of $ E_i $ and $ E_j $) uniquely as the
  intersection point $ h(E_i) \cap h(E_j) $. (It may happen that no such
  intersection exists, in which case there is no marked tropical curve of the
  given type through the points $ p_i $.)

  \begin {center} \begin{picture}(0,0)%
\includegraphics{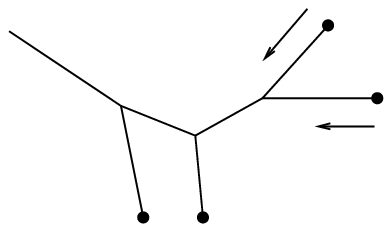}%
\end{picture}%
\setlength{\unitlength}{4144sp}%
\begingroup\makeatletter\ifx\SetFigFont\undefined
% extract first six characters in \fmtname
\def\x#1#2#3#4#5#6#7\relax{\def\x{#1#2#3#4#5#6}}%
\expandafter\x\fmtname xxxxxx\relax \def\y{splain}%
\ifx\x\y   % LaTeX or SliTeX?
\gdef\SetFigFont#1#2#3{%
  \ifnum #1<17\tiny\else \ifnum #1<20\small\else
  \ifnum #1<24\normalsize\else \ifnum #1<29\large\else
  \ifnum #1<34\Large\else \ifnum #1<41\LARGE\else
     \huge\fi\fi\fi\fi\fi\fi
  \csname #3\endcsname}%
\else
\gdef\SetFigFont#1#2#3{\begingroup
  \count@#1\relax \ifnum 25<\count@\count@25\fi
  \def\x{\endgroup\@setsize\SetFigFont{#2pt}}%
  \expandafter\x
    \csname \romannumeral\the\count@ pt\expandafter\endcsname
    \csname @\romannumeral\the\count@ pt\endcsname
  \csname #3\endcsname}%
\fi
\fi\endgroup
\begin{picture}(1784,1014)(214,-465)
\put(1759,453){\makebox(0,0)[lb]{\smash{\SetFigFont{8}{9.6}{rm}{\color[rgb]{0,0,0}$h(x_i)=p_i$}%
}}}
\put(1998,113){\makebox(0,0)[lb]{\smash{\SetFigFont{8}{9.6}{rm}{\color[rgb]{0,0,0}$h(x_j)=p_j$}%
}}}
\put(1362,115){\makebox(0,0)[rb]{\smash{\SetFigFont{8}{9.6}{rm}{\color[rgb]{0,0,0}$h(V)$}%
}}}
\put(1711,-151){\makebox(0,0)[lb]{\smash{\SetFigFont{8}{9.6}{rm}{\color[rgb]{0,0,0}$u(E_j)$}%
}}}
\put(1441,434){\makebox(0,0)[rb]{\smash{\SetFigFont{8}{9.6}{rm}{\color[rgb]{0,0,0}$u(E_i)$}%
}}}
\put(406,-241){\makebox(0,0)[b]{\smash{\SetFigFont{8}{9.6}{rm}{\color[rgb]{0,0,0}$K$}%
}}}
\end{picture}
 \end {center}

  We can now remove the two edges $ E_i $ and $ E_j $ and the marked points $
  x_i $ and $ x_j $ from $K$ and replace them by the one marked point $V$
  (whose image point $ h(V) $ is also fixed). Applying the same arguments as
  above again we conclude by induction that $ h(K) $ is uniquely defined. Of
  course, this holds for all connected components of $ |\Gamma| \backslash
  \{x_1,\dots,x_n\}$, and therefore there is at most one possibility (up to
  isomorphism) for the map $h$.
  \end {proof}

\begin {definition} \label{def-p}
  We set
    \[ \calP = \R^{2n} \backslash \bigcup_\alpha \pi(\calM_{g,\Delta}^\alpha)
       \]
  where the union is taken over all combinatorial types $ \alpha $ occurring in
  $ \overline \calM_{g,\Delta} $ such that
  \begin {itemize}
  \item $ \pi $ is not injective on $ \calM_{g,\Delta}^\alpha $; or
  \item $ \alpha $ is non-exceptional and $ \codim \alpha \ge 2 $.
  \end {itemize}
  We say that $n$ points $ p_1,\dots,p_n \in \R^2 $ are \df {in general
  position} if $ (p_1,\dots,p_n) \in \calP \subset \R^{2n} $.
\end {definition}

Our goal will be to show that the number of tropical curves (counted with a
suitable multiplicity) through $n$ given points $ p_1,\dots,p_n \in \R^2 $ in
general position does not depend on the choice of points.

\begin {remark}
  Note that our definition of points being in general position is a lot weaker
  than the definitions of \cite{Mi03} and \cite{NS04}. In \cite{Mi03} and
  \cite{NS04} points in general position require (in our language) that there
  are no curves of a combinatorial type of positive codimension through these
  points, whereas we exclude only curves of codimension at least 2. As a
  consequence, our space of points in general position will be \emph
  {connected} in $ \R^{2n} $:
\end {remark}

\begin {lemma} \label{connected} \brk
  \vspace {-1.4\baselineskip}
  \begin {enumerate}
  \item \label{connected-a}
    The locus $ \calP \subset \R^{2n} $ of points in general position is
    connected.
  \item \label{connected-b}
    The map $ \pi: \overline \calM_{g,\Delta} \to \R^{2n} $ has finite fibers
    over $ \calP \subset \R^{2n} $.
  \end {enumerate}
\end {lemma}

\begin {proof} \brk
  \vspace {-1.4\baselineskip}
  \begin{enumerate} \parindent0mm \parskip1ex plus0.3ex
  \item By definition the complement of $ \calP $ is a union of two types of
    sets:
    \begin {itemize}
    \item subspaces $ \pi(\calM_{g,\Delta}^\alpha) \subset \R^{2n} $ such that
      $ \pi $ is not injective on $ \calM_{g,\Delta}^\alpha $. As $ \pi $ is an
      affine map this means that the fiber dimension of $ \pi $ is positive.
      Moreover, by proposition \ref {linear} we must have $ \codim \alpha \ge 1
      $, and thus the dimension of $ \calM_{g,\Delta} $ is at most $ 2n-1 $
      by proposition \ref{combtype-dim}. Therefore the codimension of the image
      $ \pi(\calM_{g,\Delta}^\alpha) $ is at least 2 in $ \R^{2n} $.
    \item subspaces $ \pi(\calM_{g,\Delta}^\alpha) \subset \R^{2n} $ such that
      $ \alpha $ is non-exceptional and of codimension at least 2. Then $
      \calM_{g,\Delta}^\alpha $ has dimension at most $ 2n-2 $ by proposition
      \ref {combtype-dim}. Therefore its image has again codimension at least 2
      in $ \R^{2n} $.
    \end {itemize}
    It follows that $ \calP $ is the complement of a subset of codimension
    at least 2 in $ \R^{2n} $ and thus connected.
  \item This follows immediately since by definition $ \pi $ is injective on
    each of the (finitely many) strata $ \calM_{g,\Delta}^\alpha $ in the
    preimage of $\calP$.
  \end{enumerate}
\end {proof}

To be able to count tropical curves through points in general position we
finally need one more ingredient: we need to know the multiplicities with which
to count the curves. We will only construct these multiplicities for curves
that can actually occur through points in general position.

\begin {definition} \label{def-mult}
  Let $ (C,x_1,\dots,x_n) $ with $ C=(\Gamma,\omega,h) $ be a tropical curve
  in $ \overline \calM_{g,\Delta} $. Assume that its combinatorial type $
  \alpha $ is either exceptional or has codimension at most 1. Note that then
  all vertices of $ \Gamma $ are of valence 3 or 4.

  For a 3-valent vertex $ V \in \Gamma^0 $ with adjacent flags $ F_1,F_2,F_3
  \in \Gamma' $ we set
    \[ \mult_C V := |\det (v(F_1),v(F_2))| \in \N_{>0} \]
  (note that by the balancing condition $ v(F_1)+v(F_2)+v(F_3)= 0 $ this
  definition does not depend on how the three flags are numbered). For a
  4-valent vertex $ V \in \Gamma^0 $ with adjacent flags $ F_1,F_2,F_3,F_4 \in
  \Gamma' $ we set
    \[ \mult_C V := \max |\det (v(F_i),v(F_j)) \cdot \det (v(F_k),v(F_l))|
         \in \N_{>0} \]
  where the maximum is taken over all $ i,j,k,l $ such that $ \{ i,j,k,l \} =
  \{ 1,2,3,4 \} $. Finally, we define the \df {multiplicity} $ \mult C $ of the
  curve $C$ to be the product of the multiplicities of all its vertices.

  Using this multiplicity we can now define the function
  \begin {align*}
    N_{g,\Delta}: \qquad \qquad \quad \calP &\to \N \\
      (p_1,\dots,p_n) &\mapsto \sum \mult C
  \end {align*}
  where the sum is taken over all $ (C,x_1,\dots,x_n) \in \pi^{-1}
  (p_1,\dots,p_n) $. Note that this is a finite sum by lemma \ref{connected}
  \ref{connected-b}, and that the multiplicities $ \mult C $ are defined
  since all curves in the preimage $ \pi^{-1} (\calP) $ are by definition
  \ref{def-p} of a combinatorial type that is exceptional or of codimension at
  most 1.
\end {definition}

\begin {remark}
  In contrast to the case of counting complex curves this multiplicity with
  which we have to count curves through points in general position will \emph
  {not} always be 1 (not even with the stronger definition of general position
  of \cite{Mi03} and \cite{NS04}). It is always a positive integer though.
\end {remark}

We are now ready to prove our main theorem:

\begin {theorem} \label{constant}
  The function $ N_{g,\Delta}: \calP \to \N $ is constant. In other words, the
  number of tropical curves of genus $g$ and degree $ \Delta $ through $n$
  points in general position does not depend on the position of the points.
\end {theorem}

\begin {proof}
  As $\calP$ is connected by lemma \ref{connected} \ref{connected-a} it
  suffices to show that $ N_{g,\Delta} $ is locally constant. So let $ W
  \subset \R^{2n} $ be a small open subset of a point $ (p_1,\ldots,p_n) \in
  \calP $. Recall that $ \pi^{-1}(p_1,\dots,p_n) $ is a finite set by
  lemma \ref{connected} \ref{connected-b}. Hence as $ \pi $ is continuous and
  affine on each of the finitely many closed strata $ \overline \calM_{g,
  \Delta}^\alpha $ we can pick $W$ small enough so that the inverse image $
  \pi^{-1}(W) $ is a finite disjoint union of open subsets of $
  \overline \calM_{g,\Delta} $ each of which contains exactly one point
  in the inverse image $ \pi^{-1}(p_1,\dots,p_n) $. Of course it suffices to
  show that $ N_{g,\Delta} $ is constant when restricted to any of
  these open subsets.

  \begin {center} \begin{picture}(0,0)%
\includegraphics{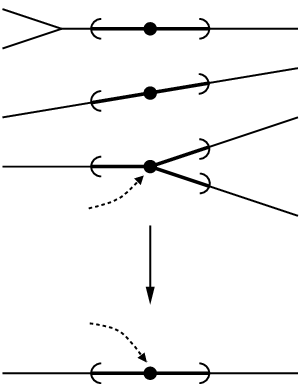}%
\end{picture}%
\setlength{\unitlength}{4144sp}%
\begingroup\makeatletter\ifx\SetFigFont\undefined
% extract first six characters in \fmtname
\def\x#1#2#3#4#5#6#7\relax{\def\x{#1#2#3#4#5#6}}%
\expandafter\x\fmtname xxxxxx\relax \def\y{splain}%
\ifx\x\y   % LaTeX or SliTeX?
\gdef\SetFigFont#1#2#3{%
  \ifnum #1<17\tiny\else \ifnum #1<20\small\else
  \ifnum #1<24\normalsize\else \ifnum #1<29\large\else
  \ifnum #1<34\Large\else \ifnum #1<41\LARGE\else
     \huge\fi\fi\fi\fi\fi\fi
  \csname #3\endcsname}%
\else
\gdef\SetFigFont#1#2#3{\begingroup
  \count@#1\relax \ifnum 25<\count@\count@25\fi
  \def\x{\endgroup\@setsize\SetFigFont{#2pt}}%
  \expandafter\x
    \csname \romannumeral\the\count@ pt\expandafter\endcsname
    \csname @\romannumeral\the\count@ pt\endcsname
  \csname #3\endcsname}%
\fi
\fi\endgroup
\begin{picture}(1452,1753)(439,-1037)
\put(1216,-511){\makebox(0,0)[lb]{\smash{\SetFigFont{8}{9.6}{rm}{\color[rgb]{0,0,0}$\pi$}%
}}}
\put(1891,-1006){\makebox(0,0)[lb]{\smash{\SetFigFont{8}{9.6}{rm}{\color[rgb]{0,0,0}$ \R^{2n} $}%
}}}
\put(1306,-871){\makebox(0,0)[b]{\smash{\SetFigFont{8}{9.6}{rm}{\color[rgb]{0,0,0}$W$}%
}}}
\put(1891,209){\makebox(0,0)[lb]{\smash{\SetFigFont{8}{9.6}{rm}{\color[rgb]{0,0,0}$\overline \calM_{g,\Delta}$}%
}}}
\put(1256,-201){\makebox(0,0)[b]{\smash{\SetFigFont{8}{9.6}{rm}{\color[rgb]{0,0,0}$U$}%
}}}
\put(816,-256){\makebox(0,0)[rb]{\smash{\SetFigFont{8}{9.6}{rm}{\color[rgb]{0,0,0}$(C,x_1,\dots,x_n)$}%
}}}
\put(816,-776){\makebox(0,0)[rb]{\smash{\SetFigFont{8}{9.6}{rm}{\color[rgb]{0,0,0}$ (p_1,\dots,p_n) $}%
}}}
\end{picture}
 \end {center}

  So let $ (C,x_1,\dots,x_n) \in \pi^{-1}(p_1,\dots,p_n) $ be a marked tropical
  curve, and let $ U \subset \pi^{-1}(W) $ be a small open neighborhood of this
  curve as above. If $ \alpha $ is the combinatorial type of $
  (C,x_1,\dots,x_n) $ then $ \alpha $ is exceptional or $ \codim \alpha \le 1 $
  by definition \ref{def-p}.

  If $ \codim \alpha = 0 $ then $ (C,x_1,\dots,x_n) $ cannot be in the
  boundary of any other stratum $ \calM_{g,\Delta}^{\alpha'} $ by proposition
  \ref{boundary} \ref{boundary-b}. Hence (after possibly shrinking $W$ and $U$)
  $ \alpha $ is the only combinatorial type occurring in $U$. As $ \dim
  \calM_{g,\Delta}^\alpha =2n $ by proposition \ref{combtype-dim} and $ \pi $
  is affine and injective on this stratum by definition \ref{def-p} we conclude
  that $ \pi|_U: U \to W $ is an isomorphism. So in this case $ N_{g,\Delta} $
  is trivially constant when restricted to $U$. In the picture above this is
  the case for the upper two connected components of $ \pi^{-1}(W) $.

  We can therefore assume from now on that $ \alpha $ is exceptional or $
  \codim \alpha = 1 $. Then $ \dim \calM_{g,\Delta}^\alpha =2n-1 $ by
  proposition \ref{combtype-dim}, and hence the image $ \pi(U \cap
  \calM_{g,\Delta}^\alpha) $ is of codimension 1 in $W$ by definition
  \ref{def-p}. We can therefore think of the image of this stratum as a
  ``wall'' that divides $W$ into two halves. We have to show that the
  restriction of $ N_{g,\Delta} $ to $U$ stays constant when crossing this
  wall. So we have to analyze which combinatorial types occur in $U$ on both
  sides of the wall. After possibly shrinking $W$ and $U$ these are of
  course just the types that contain $ \calM_{g,\Delta}^\alpha $ in their
  boundary. In the picture above this is the case for the bottom connected
  component of $ \pi^{-1}(W) $ (where we have one combinatorial type on the
  left and two combinatorial types on the right side of the wall).

  As $ \alpha $ is exceptional or $ \codim \alpha = 1 $ there are by definition
  \ref{def-combtype} four cases to check:
  \begin{enumerate} \parindent0mm \parskip1ex plus0.3ex
  \item \label{constant-a}
    $ \alpha $ is non-exceptional, and the graph $\Gamma$ has one vertex $V$
    of valence $4$:

    As this is the most interesting case (that we have already mentioned in the
    introduction) we will discuss it in detail. In fact this is the only case
    in which the number of tropical curves through the given points would not
    be constant if we did not count them with their correct multiplicities.

    Denote by $ F_1,\dots,F_4 $ the four flags with $ \partial F_i = V $ for $
    i=1,\dots,4 $, and let $ v_i := v(F_i) $ (using definition
    \ref{def-tropcurve} \ref{def-tropcurve-c}). There are exactly three
    different types $\alpha_1,\alpha_2,\alpha_3 $ that have $ \alpha $ in their
    boundary. Each of them replaces the 4-valent vertex by two 3-valent ones
    that are separated by a new edge. The remaining two edges of each of the
    two new vertices are given by the four flags $ F_1,\dots,F_4 $ in any
    possible way: we obtain type $ \alpha_i $ for $ i=1,2,3 $ when $ F_i $ and
    $ F_4 $ come together at one vertex and the other two flags at the other
    vertex. In each case the weight and direction of the new edge is determined
    uniquely by the balancing condition of definition \ref{def-tropcurve}
    \ref{def-tropcurve-c}.

    \begin {center} \input {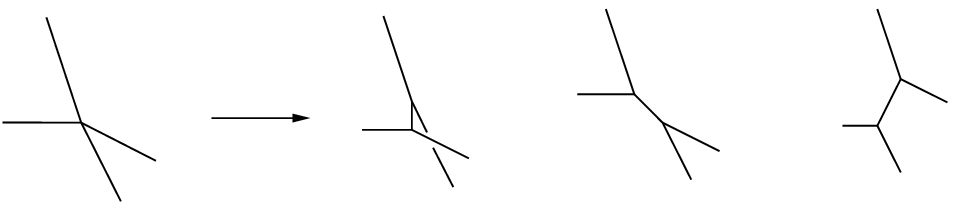} \end {center}

    Now let $ (p'_1,\dots,p'_n) \in W \backslash \pi (\calM_{g,\Delta}^\alpha)
    $ be a collection of points that is not on the wall, and let $ (C',
    x'_1,\dots,x'_n) \in \pi^{-1} (p'_1,\dots,p'_n) $. We have to determine
    which of the combinatorial types $ \alpha_1, \alpha_2, \alpha_3 $ are
    possible for $ (C',x'_1,\dots,x'_n) $. Note that in any case the new edge
    of $C'$ cannot contain a marked point. By remark \ref{rem-combtype} the
    connected component of $ |\Gamma'| \backslash \{ x'_1,\dots,x'_n\} $ that
    contains the new edge has no loops and exactly one unbounded end. Therefore
    the new edge is connected to an unbounded end via exactly one of the four
    flags $ F_1,\dots,F_4 $. By symmetry we may assume that this flag is $ F_4
    $. Then by the argument of the proof of proposition \ref{linear} for each
    possible combinatorial type the lines in $ \R^2 $ on which $ F_1 $, $ F_2
    $, and $ F_3 $ lie are fixed by the points $ p'_1,\dots,p'_n $, whereas the
    line for $ F_4 $ is not.

    Let us now determine if the collection of points $ (p'_1,\dots,p'_n) $
    admits a tropical curve of combinatorial type $ \alpha_3 $ through them.
    Let $ V = \partial F_1 = \partial F_2 $ be the common vertex of $ F_1 $ and
    $ F_2 $. As the result is obviously invariant under a relabeling of flags $
    F_1 \leftrightarrow F_2 $ we may assume without loss of generality that
      \[ \det (v_1,v_2)>0, \tag {1} \]
    i.e.\ that the (oriented) angle between the vectors $ v_1=v(F_1) $ and $
    v_2=v(F_2) $ is less than $ \pi $ (the case $ \det (v_1,v_2) = 0 $ is
    impossible since then the flags around $V$ would not span $ \R^2 $ in
    contradiction to definition \ref {def-tropcurve} \ref {def-tropcurve-c}).
    The new internal edge of the curve starts at $ h(V) $ and points in the
    direction $ -v_1-v_2 $. So a curve of type $ \alpha_3 $ through the given
    points exists if and only if the line on which $ F_3 $ lies intersects the
    ray with direction $ -v_1-v_2 $ starting at $V$. As this condition is
    invariant under a change of direction $ v_3 \leftrightarrow -v_3 $ we may
    also assume without loss of generality that
      \[ \det (v_3,v_1+v_2) > 0 \tag {2} \]
    (again this determinant cannot be zero).

    Let us consider the triangle $T$ that is cut out by the three lines on
    which the three flags $ F_1 $, $ F_2 $, and $ F_3 $ lie. The following six
    cases can occur:

    \begin {center} \input {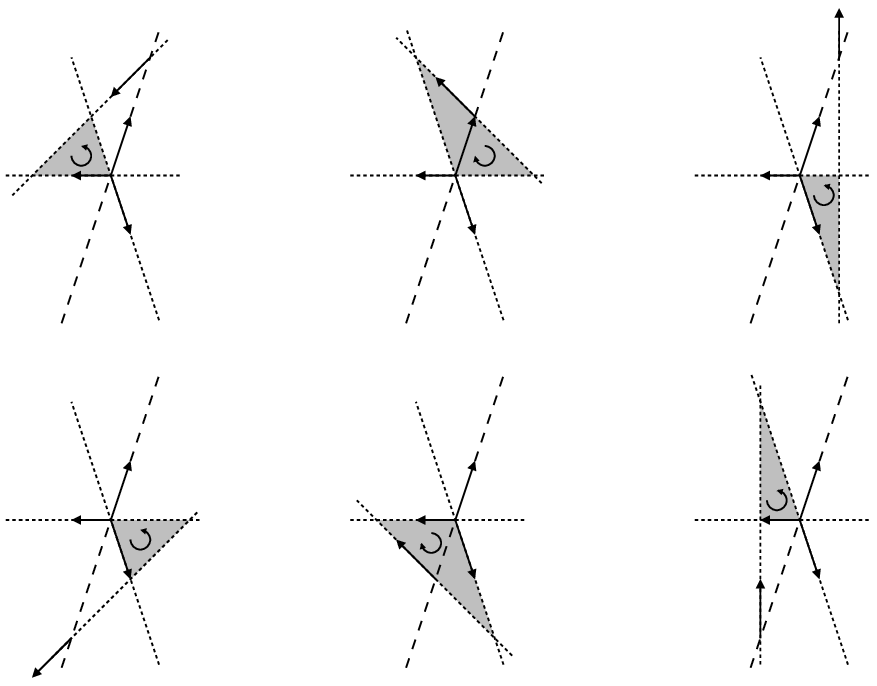} \end {center}

    In the top three cases there is a tropical curve of type $ \alpha_3 $
    through the marked points (with $ v_4 $ determined by the balancing
    condition), whereas in the bottom three cases there is no such curve. The
    cases in the left, middle, and right column differ by the slope of $ v_3 $
    compared to the slopes of $ v_1 $ and $ v_2 $.

    In any case the lines on which $ F_1 $, $ F_2 $, and $ F_3 $ lie, taken in
    this order, define an orientation on the triangle $T$ that we have also
    indicated in the picture above. Of course this also defines an orientation
    on the edges of $T$. Let $ k \in \{0,1,2,3\} $ be the number of vectors $
    v_1 $, $ v_2 $, $ v_3 $ that point against this orientation on their
    respective edge, and set
      \[ \mu := (-1)^k \cdot \prod_{\substack {V' \in \Gamma^0\\V' \neq V}}
         \mult_C V' \in \Z. \]
    By looking at all the cases above we see that there is a tropical curve of
    type $ \alpha_3 $ through the given points $ p'_1,\dots,p'_n $ if and only
    if $ \mu > 0 $ (note that the product over $ \mult_C V' $  in the
    definition of $ \mu $ is always positive by definition). To be precise it
    may also happen that the line on which $ F_3 $ lies is parallel to one of
    the other two lines. In this case the three lines do not determine a unique
    triangle but rather two ``unbounded'' triangles. One can check immediately
    that then both unbounded triangles give the same value of $ \mu $, and that
    again there is a tropical curve of type $ \alpha_3 $ through the given
    points if and only if $ \mu>0 $.

    Note that the number $ \mu $ changes sign if we exchange the labeling of $
    F_1 $ and $ F_2 $ or if we replace $ v_3 $ by $ -v_3 $. So we can remove
    our assumptions (1) and (2) above and conclude that in any case there is a
    curve of type $ \alpha_3 $ through the given points if and only if the
    number
      \[ \mu_3 := \mu \cdot \det (v_1,v_2) \cdot \det (v_3,v_1+v_2) \]
    is positive. In this case the curve is then unique by proposition
    \ref{linear}, and by definition its multiplicity is just $ \mu_3 $.

    As the number $ \mu $ is invariant under cyclic permutations $ F_1 \to F_2
    \to F_3 \to F_1 $ it now follows that for $ i=1,2,3 $ there is a curve of
    type $ \alpha_i $ through the given points if and only if $ \mu_i > 0 $,
    where
    \begin {align*}
      \mu_1 &:= \mu \cdot \det (v_2,v_3) \cdot \det (v_1,v_2+v_3),
        \quad \mbox {and} \\
      \mu_2 &:= \mu \cdot \det (v_3,v_1) \cdot \det (v_2,v_3+v_1),
    \end {align*}
    and that in this case it counts with multiplicity $ \mu_i $.

    Finally, an elementary calculation shows that $ \mu_1+\mu_2+\mu_3=0 $. From
    this it follows immediately that the number of curves in $ \pi^{-1}(p'_1,
    \dots,p'_n) \cap U $, counted with their respective multiplicities, is
      \[ \sum_{i: \mu_i>0} \mu_i = \max \{|\mu_1|,|\mu_2|,|\mu_3|\}, \]
    which is by definition the multiplicity of the one curve in $
    \pi^{-1}(p_1,\dots,p_n) $ with a 4-valent vertex.
  \item \label{constant-b}
    $ \alpha $ is non-exceptional, and the genus of $ \Gamma $ is $ g(C) = g-1
    $:

    This is a contradiction to $ \codim \alpha = 1 $ by remark \ref{genus-2}.
    Hence this case cannot occur.

  \item \label{constant-c}
    $ \alpha $ is non-exceptional, and there is a marked point $x_i$ on a
    vertex $V$ of $\Gamma$ (i.e.\ $ s(i) \in \Gamma^0 $ for one $i$):

    The idea here is the same as in case \ref{constant-a}, but the analysis is
    a lot simpler. The point $ x_i $ must lie on a 3-valent vertex $V$ whose
    adjacent flags we denote by $ F_1 $, $ F_2 $, and $ F_3 $. There are at
    most three combinatorial types that have $ \alpha $ in their boundary,
    namely the types $ \alpha_k $ for $ k=1,2,3 $ where the marked point $ x'_i
    $ lies on the edge $ [F_k] $ and the remaining data of the curve stay the
    same:

    \begin {center} \input {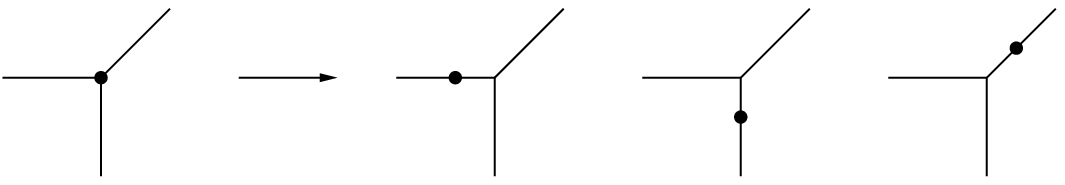} \end {center}

    By definition \ref{def-marked} \ref{def-marked-b} at least one of these
    cases will be allowed, i.e.\ lead to a marked tropical curve (namely if $
    x_i $ is moved onto the interior of the flag chosen in this definition).
    After possibly relabeling the flags we can assume that type $ \alpha_1 $ is
    possible. By remark \ref{rem-combtype} one can then reach an unbounded end
    of $ |\Gamma| \backslash \{x'_1,\dots,x'_n\} $ in $ \alpha_1 $ from both
    sides of the point $ x'_i $ along exactly one path. So one can reach such
    an unbounded end via $ F_1 $ (starting to the left of the marked point),
    and (after possibly relabeling the flags $ F_2 \leftrightarrow F_3 $) via $
    F_2 $ starting from $V$, but not via $ F_3 $ starting from $V$.

    First of all this means that type $ \alpha_3 $ is impossible, since this
    would connect the two unbounded ends behind $ F_1 $ and $ F_2 $ in
    contradiction to definition \ref{def-marked} \ref{def-marked-b}. Moreover,
    we see that (analogously to case \ref{constant-a} above) for both $
    \alpha_1 $ and $ \alpha_2 $ the line on which $ F_3 $ lies is fixed by the
    other marked points, whereas the lines on which $ F_1 $ and $ F_2 $ lie are
    not. It is now immediate that we always have exactly one of the cases $
    \alpha_1 $ and $ \alpha_2 $, depending on whether the point $ p'_i $ lies
    on the one or the other side of the line on which $ F_3 $ lies. The
    multiplicity is obviously the same in both cases as it does not depend on
    the position of the marked points.
  \item \label{constant-d}
    $ \alpha $ is exceptional:

    Note first that the two edges $ E_1 $ and $ E_2 $ joining the two 4-valent
    vertices $V$ and $V'$ are distinguishable since by definition
    \ref{def-marked} \ref{def-marked-b} at least one of them must have a marked
    point on it. At both vertices $V$ and $V'$ the balancing condition implies
    that the two other flags must point to different sides of the line on which
    $ E_1 $ and $ E_2 $ lie. We denote by $ F_1 $ and $ F'_1 $ (resp.\ $ F_2 $
    and $ F'_2 $) the flags pointing to the one (resp.\ the other) side. Then
    there are exactly two combinatorial types $ \alpha_1 $, $ \alpha_2 $ that
    contain $ \alpha $ in their boundary:

    \begin {center} \input {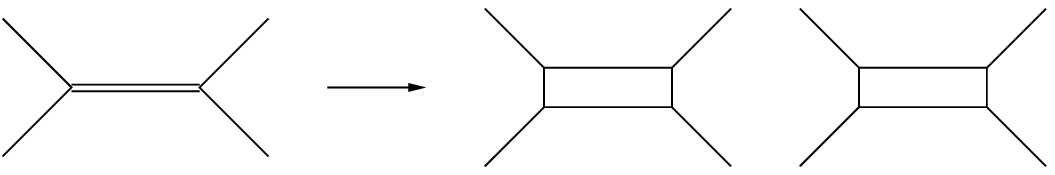} \end {center}

    Obviously, the curves of types $ \alpha $, $ \alpha_1 $, and $ \alpha_2 $
    all have the same multiplicity. Applying the analysis of \ref{constant-a}
    to both vertices $V$ and $V'$ we see that the types $ \alpha_1 $ and $
    \alpha_2 $ always occur on different sides of the wall, whereas $ \alpha $
    occurs on the wall itself. Hence it follows that restriction of $
    N_{g,\Delta} $ to $U$ is constant.
  \end{enumerate}
\end {proof}

\begin{remark} \label{mik}
  Let us finally recap what we have shown and how our results relate to the
  existing literature.

  Let $X$ be a (complex) smooth projective toric surface defined by a fan $
  \Sigma $, and let $ G' $ be the sub-semigroup of $G$ (as introduced in
  definition \ref{def-degree}) generated by all primitive integral vectors
  along the rays of $ \Sigma $. It is then well-known that every element $
  \Delta = u_1 \oplus \cdots \oplus u_k \in G' $ with $ u_1+\cdots+u_k=0 $
  corresponds to a homology class of complex curves in $X$. Pick a genus $ g
  \ge 0 $, and set $ n = \#\Delta+g-1 $. Moreover, choose $n$ points $
  q_1,\dots,q_n \in X $ in general position that lie in the torus $ (\C^*)^2
  \subset X $. Define the map $ \Log: (\C^*)^2 \to \R^2 $ by $ \Log(z_1,z_2)
  := (\log |z_1|,\log |z_2|) $ and set $ p_i := \Log (q_i) $ for $ i=1,\dots,n
  $. It is then the main result of \cite{Mi03} (the so-called ``Correspondence
  Theorem'') that the number of \emph{complex} curves in $X$ of class $ \Delta
  $ and genus $g$ through the points $ q_1,\dots,q_n $ is equal to the number
  of \emph{tropical} plane curves of degree $ \Delta $ and genus $g$ through
  the points $ p_1,\dots,p_n $ (counted with their multiplicities as
  constructed in definition \ref{def-mult}). In particular, since the number of
  such complex curves does not depend on the points $ q_1,\dots,q_n $, it
  follows as a corollary that the number of such tropical curves does not
  depend on the position of the points $ p_1,\dots,p_n $ either. In this paper
  we have given a proof of this last statement within tropical geometry, i.e.\
  without referring to the (highly non-trivial) Correspondence Theorem.

  For example, if we choose for $ \Delta $ the degree that contains $d$ times
  the vectors $ (-1,0) $, $ (0,-1) $, and $ (1,1) $ each then tropical curves
  of degree $ \Delta $ correspond to complex curves in $ \P^2 $ of degree $d$.
\end{remark}

\begin {remark} \label{enden}
  Our results are in fact more general than what we have just described in
  remark \ref{mik}. Namely, our main theorem \ref{constant} is also applicable
  in the following two cases in which there is no analogue of the
  ``Correspondence Theorem'' yet:
  \begin {enumerate}
  \item \label{enden-a}
    if the degree $ \Delta $ does not only contain primitive integral
    vectors, i.e.\ if there are unbounded ends with weights greater than 1;
  \item \label{enden-b}
    if some of the unbounded ends are ``fixed'', i.e.\ if not only their slope
    is fixed but also the line in $ \R^2 $ on which they lie. We have not
    included this set-up explicitly in our definitions since this would have
    made the notations too complicated. The necessary modifications in the
    constructions and proofs are very straightforward however if one thinks of
    a fixed unbounded end as an unbounded end with a marked point ``at
    infinity'' on it. For example, every fixed unbounded end reduces the
    required number $n$ of marked points by 1, and in definition
    \ref{def-marked} \ref{def-marked-b} we have to replace ``no connected
    component with more than one unbounded end'' by ``no connected component
    with more than one unbounded end that is not fixed''.
  \end {enumerate}
  Intuitively, case \ref{enden-a} corresponds to curves with fixed
  multiplicities to the corresponding toric divisors (e.g.\ tropical curves of
  degree $ (-2,0) \oplus (0,-1) \oplus (0,-1) \oplus (1,1) \oplus (1,1) $
  correspond to complex conics in $ \P^2 $ tangent to the line determined by
  the vector $ (-1,0) $). Case \ref{enden-b} corresponds to fixed intersection
  points of the complex curves with these toric divisors. Combining both cases
  we should in some cases be able to construct a tropical analogue of relative
  Gromov-Witten invariants, i.e.\ of complex curves with fixed multiplicities
  to a given divisor in maybe fixed points of this divisor. We will discuss
  this in detail in a forthcoming paper \cite{GM05}.
\end {remark}

\bibliographystyle{amsalpha}
 \bibliography{bibliographie}

\end{document}